%% file: sl3rev.tex
\documentclass{amsart}
\usepackage{amssymb}
\usepackage{amscd}
\usepackage{verbatim}
\usepackage{epsfig}
\newcommand{\bM}{\overline{M}}
\newcommand{\res}{{\mathrm{res}}}
\newcommand\calH{\mathcal{H}}
\newcommand\moddist{\tilde\delta}
\begin{document}

\input macro.tex

\renewcommand{\theenumi}{\roman{enumi}}
\renewcommand{\labelenumi}{(\theenumi)}

\title[Scattering theory on $\SL(3)/\SO(3)$]
{Scattering theory on $\SL(3)/\SO(3)$: \\ 
connections with quantum $3$-body scattering}
\author[Rafe Mazzeo and Andras Vasy]{Rafe Mazzeo and Andr\'as Vasy}
\address{R.\ M.: Department of Mathematics, Stanford University, Stanford,
CA 94305}
\email{mazzeo@math.stanford.edu}
\address{A.\ V.: Department of Mathematics, Massachusetts Institute of
Technology, MA 02139}
\email{andras@math.mit.edu}
\date{June 10, 2002. Revised November 26. 2002.}

\begin{abstract}
In this paper we continue our program of extending the methods of 
geometric scattering theory to encompass the analysis of the Laplacian
on symmetric spaces of rank greater than one and their geometric
perturbations. Our goal here is to explain how analysis of the Laplacian 
on the globally symmetric space $\SL(3,\RR)/\SO(3,\RR)$ is very closed 
related to quantum three-body scattering. In particular, we adapt geometric 
constructions from recent advances in that field by one of us (A.V.), 
as well as from our previous work \cite{Mazzeo-Vasy:Resolvents} concerning
resolvents for product spaces, to give a precise description of the
resolvent and the spherical functions on this space. Amongst the many
technical advantages, these methods give results which are uniform up to 
the walls of the Weyl chambers. 
\end{abstract}

\maketitle

\section{Introduction}
It has long been observed that there are formal similarities between
the spectral theory for Laplacians on (locally and globally)
symmetric spaces of rank greater than one and Hamiltonians associated
to quantum $N$-body interactions. Our contention is that these 
similarities have deep-seated explanations, rooted in the geometry
of certain natural compactifications of the spaces involved and the 
asymptotic structure of these operators, and 
that the methods of geometric scattering theory constitute a natural
set of techniques with which to study both problems. In the present paper 
we use these methods to provide an alternate perspective on mostly 
well-known results concerning the Laplacian on the globally symmetric space
\[
M = \SL(3,\RR)/\SO(3,\RR).
\]
Besides giving a new set of methods to study scattering theory on this space
which are not constrained by the algebraic rigidity and structure, this 
more general approach has benefits even in this classical framework.
Specifically, starting from the perturbation expansion methods of 
Harish-Chandra, as explained in 
\cite{Helgason:Groups}, and continuing through recent developments by
Anker and Ji \cite{Anker:Forme}, \cite{Anker-Ji:Comportement},
\cite{Anker-Ji:Heat}, it has always been problematic to obtain uniformity
of various analytic objects near the walls of the Weyl chambers. We obtain
this uniformity as a simple by-product of our method.

Let us now briefly set this work in perspective. The recent advances
in quantum $N$-body scattering from the point of view of geometric
scattering, to which we alluded above, are detailed in
\cite{Vasy:Structure}, \cite{Vasy:Bound-States} and \cite{Vasy:Propagation-2},
and we shall not say much more about this work here. Next, there are 
very many applications of geometric scattering theory to scattering 
on asymptotically Euclidean spaces and
locally and globally symmetric spaces of rank one, \cite{RBMSpec}, 
\cite{Hassell-Vasy}, \cite{Mazzeo-Melrose:Meromorphic}, 
\cite{Epstein-Melrose-Mendoza:Resolvent}, to name just a very few 
(and concentrating on those most relevant to the present discussion). 
More recently there has been progress on geometric scattering on higher 
rank spaces. For example, Vaillant \cite{Vaillant} has extended M\"uller's 
well-known $L^2$ index
theorem for spaces with ${\mathbb Q}$-rank one ends to a general
geometric setting. Most germane to the present work is 
\cite{Mazzeo-Vasy:Resolvents}, which contains the beginnings of 
a serious approach to dealing with the main technical problem of `corners 
at infinity' which arise in higher rank geometry. That paper focuses
on the special cases of products of hyperbolic, or more generally,
asymptotically hyperbolic spaces, and produces a thorough analysis 
of the resolvent of the Laplacian on such spaces, including such
features as its meromorphic continuation and the fine structure of
its asymptotics at infinity. This analysis includes the construction
of a geometric compactification of the double-space of the product
space, on which the resolvent naturally lives as a particularly
simple distribution, and which we call the resolvent double space.
The methods of that paper rely heavily on the product structure,
and an interesting representation formula for the resolvent
in terms of the resolvents on the factors which is afforded by this 
structure. While not perhaps apparent there, the final results are 
in fact independent of this product structure and obtain in much
more general situations.  

Before embarking on a general development of the analysis of the 
resolvent for spaces with `asymptotically rank two (or higher) geometry', 
we have thought it worthwhile to explain in detail how these methods apply 
to this specific example, $M=\SL(3,\RR)/SO(3,\RR)$, since it is a natural 
model space for the more general situation and is of substantial 
interest in and of itself. The methods here should apply more generally
without any new ideas, just a bit more sweat and tears! Our aim is several-fold. 
At the very least we wish to emphasize the resolvent double-space, which
is a compactification of $M \times M$ as a manifold with corners, 
and its utility for obtaining and most naturally phrasing results about 
the asymptotics of the 
resolvent; we also wish to show how the seemingly special product analysis 
of \cite{Mazzeo-Vasy:Resolvents} emerges as the `model analysis' in this 
non-product setting. 

Let us now describe our results in more detail. Fix an invariant 
metric $g$ on $M$, and let $\Delta_g$ and $R(\ev) = (\Delta_g-\ev)^{-1}$
denote its Laplacian and resolvent. We wish to examine the structure
of the Schwartz kernel of $R(\ev)$, and in particular to determine its
asymptotics as the spatial variables tend to infinity in $M$. We do this 
here for $\ev$ in the resolvent set, and with additional work (which is
not done here) this can also be carried out for $\ev$ approaching 
$\mbox{spec}\,(\Delta_g)$. As in the 
traditional approach, the invariance 
properties of $\Delta$ allow us to reduce the analysis to that on the
flats, and this turns out to be very close to three-body scattering.
The central new feature of the analysis is to replace the perturbation series
expansion of Harish-Chandra by $L^2$-based scattering theory methods 
in the spirit of \cite{Vasy:Structure}. As noted earlier, the results
are easily seen to be uniform across the walls of the Weyl chambers. 
In future work we shall study the resolvent for spaces with `asymptotically
rank-two geometry', which is only slightly more difficult; unlike there,
however, the present analysis is an explicit mixture of algebra 
(the reduction) and geometric scattering theory (the three-body problem). 

In order to describe the structure of $R(\ev)$, we first define a
compactification $\olM$ of $M$ itself. Recall that $M$ is identified with 
the set of $3$-by-$3$ positive definite matrices of determinant 1; it
is five-dimensional, and its compactification $\olM$ is a $\Cinf$ manifold 
with corners of codimension two. $\olM$ has two boundary hypersurfaces,
$H_\sharp$ and $H^\sharp$, in the interior of each of which the ratio 
of the smaller two, respectively the larger two, eigenvalues of the
representing matrix is bounded. Correspondingly, either of these boundary 
faces is characterized by the fact that the ratio of the appropriate two 
eigenvalues extends to vanish on that face, hence gives a local boundary 
defining function. The subspace of diagonal matrices is identified with the 
flat $\exp(\fraka)$, and the Weyl group $W= S_3$ acts on it by 
permutations; its closure in $\olM$ is a hexagon, the faces of which
are permuted by the action of $W$. The fixed point sets of elements of 
the Weyl group partition $\fraka$ into the Weyl chambers; 
the fixed point sets themselves constitute the Weyl chamber walls, and the
closure of the chambers in $\olM$ are the sides of the hexagon. 
Adjacent sides of the hexagon lie in different boundary hypersurfaces of 
$\olM$. The boundary hypersurfaces of $\olM$ are equipped with a fibration 
with fibers $\SL(2,\RR)/\SO(2,\RR)=\HH^2$. For example, two interior points 
of $H^\sharp$ are in the same fiber if the sum of the eigenspaces of the two 
larger eigenvalues (whose ratio is, by assumption, bounded in this region) 
is the same. This gives $\olM$ a boundary fibration structure, 
similar to (but more complicated than) ones considered in 
\cite{RBMSpec, Mazzeo:Edge, Mazzeo-Melrose:Fibred}.

We must blow up $\olM$ further in order to describe the resolvent 
efficiently. The motivation for this is that the flat is simply
a Euclidean space, and its most natural compactification is therefore
the usual radial one, also known as the geodesic compactification. 
We wish to find a compactification of
$\olM$ compatible with this. Let $\rho_\sharp$ and $\rho^\sharp$ denote 
boundary defining functions for $H_\sharp$ and $H^\sharp$ (defined
in terms of ratios of eigenvalues, as above). 
Replace these with the `slow variables' $-1/\log\rho_\sharp$ and 
$-1/\log\rho^\sharp$, respectively; by definition, $\olM$ with
these new boundary defining functions is denoted $\olM_{\log}$, the 
logarithmic blow-up of $\olM$. Now perform the spherical blow-up of 
the corner $H_\sharp \cap H^\sharp$ in this space. This sequence of 
operations results in the final `single space'
\[
{\widetilde M}=[\olM_{\log};H_\sharp\cap H^\sharp]
\]
To check that this meets the requirement of compatibility
with the radial compactification of the flat, note that if
$r$ is a Euclidean radial variable on $\fraka$ (outside a compact set), 
then its inverse $r^{-1} \equiv x\in\Cinf(\widetilde M)$ 
is a total boundary defining function of $\widetilde M$, i.e.\ 
vanishes simply on all faces. 

\begin{figure}[ht]
\begin{center}
\mbox{\epsfig{file=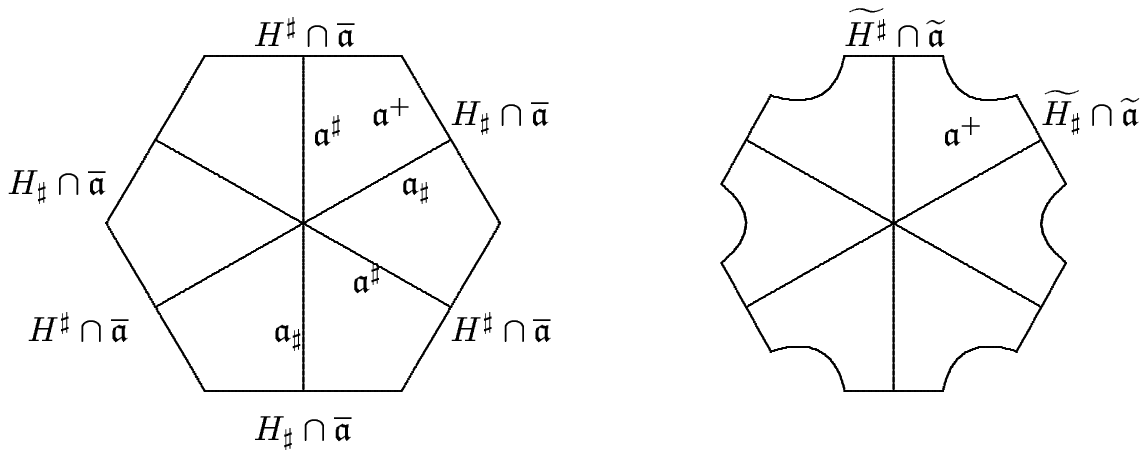}}
\end{center}
\caption{The closure of $\fraka$ in the compactifications
$\bar M$ and $\tilde M$ of $M$.}
\label{fig:flat2}
\end{figure}

A standard preliminary result in scattering theory concerns the
far-field behaviour of $R(\ev)f$ where $f \in \Cinf_c$, initially
when $\ev$ lies in the resolvent set, and later when it approaches
the spectrum, cf.\ \cite{RBMGeo}. The geometry of $\tilde{M}$ has
been set up precisely so that the analogous result here appears
quite simple. 
\begin{thm*}
Suppose $f\in\Cinf_c(X)$, and $\ev\nin\spec(\Delta)$.
Then, with $\lambda_0 = 1/3$, 
\[
R(\ev)f=\rho_\sharp\rho^\sharp x^{1/2} \exp\left(-i\sqrt{\ev-\ev_0}/x\right)
g, 
\]
where $g\in\Cinf({\widetilde M})$ vanishes to first order on the lifts of 
$H_\sharp$ and $H^\sharp$ to $\tilde M$, but is nonvanishing on the other
hypersurface (obtained by blowing up $H_\sharp \cap H^\sharp$). 
The square root in the exponential is the one having negative imaginary 
part in the resolvent set $\lambda \in \CC \setminus [\lambda_0,\infty)$. 
If $f\in\dist_c(X)$, then a similar statement holds (away from $\ssupp f$). 
\end{thm*}

\begin{rem} As already indicated, for $\ev \to \text{spec}\,(\Delta)$
this result extends as
$$
R(\ev-i0)f=\rho_\sharp\rho^\sharp x^{1/2} e^{-i\sqrt{\ev-\ev_0}/x} g,
\qquad g\in\Cinf(\tilde M).
$$
The basic strategy of the proof of this extension follows the same
lines, but exactly as in the three-body setting it requires a somewhat
more elaborate phase space analysis. For simplicity we do not discuss this 
extension in the present paper.
\end{rem}

The leading coefficients in the expansions of $R(\ev)f$ at the various
boundary faces, i.e.\ the values of $g$ or its first derivative, may
be calculated explicitly, and this leads in a straightforward manner
to a description of the Martin compactification of $\SL(3,\RR)/\SO(3,\RR)$
as $\tilde M$, with its smooth structure removed. This identification
of the Martin boundary is due originally to Guivarc'h, Ji and Taylor
\cite{Guivarch-Ji-Taylor:Compactifications}, using estimates of Anker and Ji
\cite{Anker:Forme, Anker-Ji:Comportement} to control the (previously 
opaque) behavior of the resolvent kernel at the Weyl chamber
walls. In fact, the estimates of \cite{Anker-Ji:Comportement}, 
see also \cite[Section~8.10]{Guivarch-Ji-Taylor:Compactifications},
amount to upper and lower bounds for $R(\ev)f$, when $\ev$ is real and in
the resolvent set, by expressions of the same form as in our theorem. 
In later work, Anker and Ji describe the leading term of the asymptotics
uniformly by algebraic methods. Our method is more analytic, and since
it automatically gives uniform asymptotics, the Martin compactification 
can be directly read off from it, just as in our previous work 
\cite{Mazzeo-Vasy:Resolvents}. 

There is another approach, due to Trombi and Varadarajan 
\cite{Trombi-Varadarajan:Spherical}, which is intermediate between
our approach and that of Harish-Chandra. In their approach one
constructs spherical functions as sums of polyhomogeneous conormal 
functions on $\olM$ by constructing their Taylor series at all boundary
hypersurfaces of $\olM$. By comparison, Harish-Chandra's method
amounts to constructing the spherical functions in Taylor series at the
corner of $\olM$. Owing to the algebraic nature of the space, 
these Taylor series actually converge in the appropriate 
regions; this does not hold in more geometric settings, of course.

As another application of our resolvent estimates, we also take
up the construction of the spherical functions. 
This construction is essentially just that of Trombi and 
Varadarajan, but instead of appealing to convergence of the Taylor 
series, we use the resolvent to remove the error term, 
and this results in an additional term with the same asymptotics as 
the Green function. However, since the Taylor series actually 
converges, the error term vanishes, and hence the Green function asymptotics
do not appear in the asymptotics of the spherical function; this is 
the extent to which algebra enters into our analysis.

To set this discussion into the language of Euclidean scattering, and
in particular to compare with the language of three-body 
scattering, the spherical functions on $M$ are analogues of (reflected) 
`plane waves' on the flats, corresponding to colliding particles,
although here the eigenvalues collide; on the other hand, the Green 
function for $\Delta$ on $M$ is the analogue of a `spherical wave' in
Euclidean scattering. The conflict of terminology is somewhat unfortunate.

\bigskip

\noindent{\bf Overview of the parametrix construction}

\noindent We now sketch in outline some details of our methods and 
constructions.  The goal is to construct the Schwartz kernel of
the resolvent as a distribution with quite explicit singular structure
on some compactification of $\olM^2$. This is accomplished by
constructing a sequence of successively finer approximations to
$(\Delta - \ev)^{-1}$, where `fineness' is measured by the extent
to which these operators map into spaces with better regularity and 
decay at infinity. These parametrices lie in certain `calculi' of
pseudodifferential operators which are defined by fixing the possible 
singular structures of the Schwartz kernels of their 
elements both at the diagonal, but more interestingly, near
the boundary of $\olM^2$. 

The first step is the construction of
a parametrix in the `small calculus', which we also call the 
edge-to-edge calculus, of pseudodifferential operators on $M$.
In fact, this is defined on any manifold with corners up to codimension 
two which has fibrations on its boundary faces analogous to those of 
$\SL(3)/\SO(3)$. A more general development of this calculus will appear 
elsewhere. The parametrix $G(\ev)$ for $\Delta-\ev$ in this 
calculus has the property that the error $E(\ev)=G(\ev)(\Delta-\ev)-\Id$ 
is smoothing but does not increase the decay rate of functions, hence is 
not compact on $L^2(M,dV_g)$. The constructions within this
small calculus are merely a systematic way of organizing the local 
elliptic parametrix construction uniformly to infinity, and this
parametrix gives scale-invariant estimates uniform to $\pa M$. To
amplify on this last statement, we may use this calculus to define
the Sobolev spaces $\Hss^m(\olM) = \{u\in L^2(M,dV_g):\ 
\Delta^{m/2} u\in L^2(M,dV_g)\}$, $m\geq 0$. These spaces 
reflect the basic scaling structure of $\olM$ near its boundaries. 

At this point we use the group structure to simplify matters by effectively 
reducing to the flat. For $p \in M$, let $K_p$ denote the subgroup of 
$\SL(3)$ fixing this point; we may as well assume that $p$ is the 
identity matrix, which identifies the subspace of $K_p$-invariant functions
with the space of Weyl group invariant functions on the flat $\exp(\fraka)$
(or equivalently, of functions on diagonal matrices invariant under 
permutation of the diagonal entries). Since the Green function
(for $\Delta - \ev$) with pole at $p$ is $K_p$-invariant, we may
as well consider only parametrices which respect this structure. 
This reduction is certainly helpful, but not essential; it is the key 
point where our restriction to the actual symmetric space makes
a difference in terms of simplifying the presentation. 

Denote by $\Hss^m(\olM)^{K_p}$ the invariant elements of the Sobolev space 
$\Hss^m(\olM)$. The initial parametrix $G(\ev)$ constructed in
the first step may not preserve $K_p$ invariance, but this is
easily remedied by averaging it over $K_p$; this produces an operator 
$G_p(\ev)$ which satisfies
\[
\tilde E(\ev)=G_p(\ev)(\Delta-\ev)-\Id:L^2(M,dV_g)^{K_p}
\longrightarrow \Hss^m(\olM)^{K_p} \qquad \text{for all}\  m.
\]

A second step is needed to get a parametrix with a decaying (as well
as smoothing) error term. Thus we wish to construct another
parametrix $\tilde R(\ev)$ for $\Delta-\ev$, $\ev\nin\spec(\Delta)$, 
which acts on these $K_p$-invariant function spaces and satisfies
\[
\tilde R(\ev)(\Delta-\ev)-\Id:\Hss^m(\olM)^{K_p}\longrightarrow 
x^s \Hss^m(\olM)^{K_p}\qquad \text{for all}\  s.
\]
Granting this for a moment, we combine these two operators
to get
\begin{equation*}
\left(G_p(\ev)+\tilde E(\ev)\tilde R(\ev)\right)(\Delta
-\ev)-\Id:L^2(M,dV_g)^{K_p}\to x^s\Hss^m(\olM)^{K_p}
\end{equation*}
for any $s > 0$. This error term is now compact on $L^2(M,dV_g)^{K_p}$,
and so using the simplest spectral properties of $\Delta$, we may 
remove it and obtain an inverse to $\Delta-\ev$ acting on $K_p$-invariant 
functions. Since, as remarked before, $(\Delta - \ev)^{-1}$ is
necessarily $K_p$-invariant, we have captured the full resolvent.

\begin{figure}[ht]
\begin{center}
\mbox{\epsfig{file=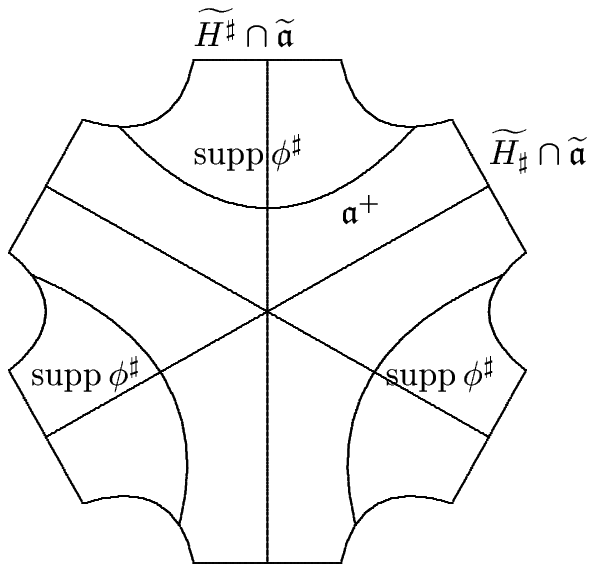}}
\end{center}
\caption{The intersection of $\supp\phi^\sharp$ with the flat
$\widetilde{\fraka}$.}
\label{fig:flat3}
\end{figure}

The main subtleties in this paper center on the construction of 
$\tilde R(\ev)$, to which we now turn. 
We return to the single space $\widetilde{M}$.
Denote by $\widetilde{H_\sharp}$ and $\widetilde{H^\sharp}$ the lifts of 
the boundary faces $H_\sharp$ and $H^\sharp$ from $\olM$ to
$\widetilde{M}$, and let $\fraka_\sharp$ and $\fraka^\sharp$ be the 
Weyl chamber walls intersecting these faces, respectively. Choose
a $\Cinf({\widetilde M})^{K_p}$ partition of unity, $\phi_\sharp+\phi^\sharp+
\phi_0=1$ on ${\widetilde M}$ such that $\supp\phi_\sharp$ is disjoint 
from $\widetilde{H^\sharp}\cap \fraka^\sharp$ and $\supp\phi^\sharp$ is 
disjoint from $\widetilde{H_\sharp}\cap \fraka_\sharp$, and with 
$\phi_0 \in \Cinf_c(M)$. (These can be constructed on the closure of
$\widetilde{\fraka}$ and extended to $K_p$-invariant functions on 
$\widetilde M$.) Let $\psi_\sharp$, $\psi^\sharp$ be $K_p$-invariant 
cutoffs which are 
identically $1$ on $\supp\phi_\sharp$ and $\supp\phi^\sharp$, respectively,
and which vanish on $\widetilde{H^\sharp}\cap \fraka^\sharp$ and 
$\widetilde{H_\sharp} \cap\fraka_\sharp$, respectively. 

Along $\widetilde{H^\sharp}$, $\Delta$ is well approximated by the 
product operator $L^\sharp=\frac{3}{4}(sD_s)^2+i\frac{3}{2}(sD_s)+
\Delta_{\HH^2}$, where $s=\rho^\sharp$ and $\Delta_{\HH^2}$ is the 
Laplacian on the fiber $\HH^2$ of $\widetilde{H^\sharp}$. More precisely, 
\[
\Delta-L^\sharp:\Hss^m(\olM)^{K_p}\longrightarrow 
\rho^\sharp\Hss^{m-1}(\olM).
\]
There is an analogous product operator $L_\sharp$ which approximates
$\Delta$ near $\widetilde{H_\sharp}$. One small complication is 
that because of its structure at $\fraka_\sharp$, $L^\sharp$ does not 
preserve $K_p$-invariance, but this is not serious since  
$\psi^\sharp L^\sharp \phi^\sharp$ does preserve this invariance.
The fact that we can approximate $\Delta$ by product operators is
the big gain, and is one of the remarkable things accomplished
by passing to $\widetilde M$, for in fact the representation formula 
for the resolvent on a product space from \cite{Mazzeo-Vasy:Resolvents} 
gives a precise description of $(L^\sharp-\ev)^{-1}$ and 
$(L_\sharp-\ev)^{-1}$. We use these as local models for
the structure of the improved parametrix, and put
\begin{equation*}
\tilde R(\ev)=\psi^\sharp(L^\sharp-\ev)^{-1}\phi^\sharp+\psi_\sharp
(L_\sharp-\ev)^{-1}\phi_\sharp.
\end{equation*}
It is not difficult to see that this has all the desired properties.

As noted earlier in this introduction, this construction
can be extended to let $\ev$ approach the spectrum, so as
to obtain the structure of the limiting values of the resolvent 
$R(\ev\pm i0)= (\Delta-(\ev\pm i0))^{-1}$. This is not done here 
for reasons of brevity only. We comment briefly, in this paragraph
and the next, on some consequences of this extension.  
The issue is that we must keep track of the `propagation of singularities'
along $\del \widetilde{\exp(\fraka)}$, as in in three-body scattering 
\cite{Vasy:Structure}. In this context, the term `singularity'
refers to a microlocal description of the lack of rapid decay at 
infinity. An explicit iteration allows us to construct successively
finer parametrices, leaving error terms which map $L^2(M,dV_g)^{K_p}\to 
x^k\Hss^m(\bar M)^{K_p}$ for higher and higher values of $k$ and $m$. 
The terms in this iterative series correspond to singularities reflecting 
from the walls at infinity. In fact, just as in three-body scattering, 
there are only finitely many (in fact three) such reflections. 

The spherical functions centered at $p$ are also $K_p$ invariant,
and are parametrized by incoming directions $\xi$, $|\xi|^2 = 
\ev-\ev_0$. These are constructed as perturbations of the plane waves 
\begin{equation*}
u_\xi(z)=\rho_\sharp(z)\rho^\sharp(z)e^{-i \xi\cdot z}
\end{equation*}
on $\fraka$ (which we identify with $\exp(\fraka)$), 
where $z$ a Euclidian variable. In fact, on $\tilde M$, 
$(\Delta-\ev)u_\xi$ decreases rapidly away from $\widetilde{H^\sharp}$
and $\widetilde{H_\sharp}$. 
More importantly, if $\xi \notin \fraka_\sharp \cap \fraka^\sharp$,
then $(\Delta-\ev)u_\xi$ is nowhere incoming, in the sense of the 
scattering wave front set, so that $R(\ev+ i0)$ can be applied to it. 
The detailed structure of $R(\ev+i0)$ discussed above leads to reflected 
plane waves. There are six such terms, corresponding 
to the six elements of the Weyl group. These correspond precisely to 
the six terms in Harish-Chandra's construction of spherical functions. 
The coefficients of the leading terms, which are Harish-Chandra's 
${\mathfrak c}$-function, correspond to the scattering matrices of the 
`two-body problems', in this case the scattering matrices on $\HH^2$.
As before, this analysis would allow us to let $\xi$ approach the walls.
Certain aspects of this still appear in the construction of
off-spectrum spherical functions in \S 6 below. 

\bigskip

The remainder of this paper is organized as follows. In \S 2 we review
the geometry of $M=\SL(3)/\SO(3)$ and its compactification $\olM$. 
The small calculus of pseudodifferential operators on $M$ is defined 
in \S 3 through the properties of the Schwartz kernels of its elements 
on the resolvent compactification of $M \times M$. This leads to
the first parametrix for $\Delta-\ev$, which captures the diagonal 
singularity of $R(\ev)$ uniformly to infinity. 
In \S 4, we discuss a model problem on $\Real\times\HH^2$, 
which is used for the construction of the finer parametrix
in \S 5. In \S 6 we consider spherical functions, and discuss
the extent to which algebra plays a role in their asymptotics. 
The Appendix contains a summary of results from 
\cite{Mazzeo-Vasy:Resolvents} concerning resolvents for product problems.

The authors are grateful to Lizhen Ji and Richard Melrose for helpful
discussions and for encouragement. They also thank the Mathematical
Sciences Research Institute where part of the work was completed during
the semester-long program in Spectral and scattering theory in Spring 2001.
R.\ M.\ is partially supported by NSF grant \#DMS-991975;
A.\ V.\ is partially supported by NSF grant \#DMS-9970607.
A.\ V.\ also thanks the Erwin Schr\"odinger Institute
for its hospitality during his stay in Vienna, Austria, while working
on this paper.

\section{Geometric preliminaries}
We refer to \cite{Eberlein:Geometry}, \cite{Helgason:Groups} 
and \cite{Jost:Riemannian} for a nice general discussion of
non-compact symmetric spaces, though all the essential ingredients 
are discussed below.

\subsection{Geometry}
Any element $A \in \SL(3)$ admits a unique polar decomposition 
$A=VR$, where $V=(A A^t)^{1/2}$ is positive definite and symmetric,
with $\det V = 1$, and $R\in \SO(3)$. This yields the standard identification 
of the symmetric space $\SL(3)/\SO(3)$ with the set $M$ of positive 
definite $3$-by-$3$ matrices of determinant $1$, where 
\[
\SL(3)/SO(3) \ni A\cdot \SO(3) \longmapsto (A A^t)^{1/2} \in M.
\]
We write $[A]$ for the image of $A$ in $M$. 
The action of $\SL(3)$ on $M$ is described by 
\[
\SL(3) \times M \ni (B,V) \longmapsto \phi_B(V) = (B V^2 B^t)^{1/2}.
\]
Thus $(\phi_B)_*$ identifies $T_{\Id}M$ with $T_{[B]}M$. 
$T_{\Id} M$ consists of symmetric matrices of trace $0$, and we
shall use the Killing form $g(A,A)=6\Tr (A A^t) = 
6\sum_{ij} a_{ij}^2$ as the metric on this subspace.
This induces a unique invariant Riemannian metric $g$ on all of $M$
in the usual way.

We need to analyze the structure of $M$ at infinity, and so we proceed
as follows. This space is almost a product, except along some lower 
dimensional strata. To understand this, diagonalize any $B \neq \Id$, i.e.\ 
write $B=O\Lambda O^t$, where $\Lambda$ is diagonal with positive
entries and $O\in\SO(3)$. We write $\fraka$ for the subspace of
diagonal matrices of trace zero, so that $\exp(\fraka) = M \cap \Diag_3$ 
is the space of diagonal matrices with positive entries. Neither $O$ nor 
$\Lambda$ are uniquely determined here since the ordering of the entries 
of $\Lambda$ is not fixed. Let $\calP$ denote the subgroup of all signed 
permutation matrices in $\SO(3)$, i.e.\ orthogonal matrices of
determinant $1$ which preserve the 
subspace decomposition $\RR^3 = \RR \oplus \RR \oplus \RR$ rather than just
permuting basis elements. If $P \in \calP$ then $P\Lambda P^{t}$ is 
again diagonal with positive entries. Conversely, $\Lambda$ and $O$ are 
determined up to the action of $P\in\calP$. Here $\calP$ acts on $\Lambda$
by conjugation and on $\SO(3)$ by right multiplication by the inverse:
$O\Lambda O^t=(OP^t)(P\Lambda P^t)(OP^t)^t$.
Let $(\exp(\fraka))^*$ be the subset of $\exp(\fraka)$
where the diagonal entries are all distinct. The action on $\SO(3)$,
hence on the product $(\exp(\fraka))^*\times\SO(3)$ is
free, so $((\exp(\fraka))^*\times\SO(3))/\calP$ may be identified
with $M\setminus \calC$, where $\calC$ is the subset of elements $[A] 
\in M$ where at least two of the eigenvalues coincide.

Let $\calP'$ be the subgroup of
$\calP$ that preserves each summand $\RR$ in
$\RR^3 = \RR \oplus \RR \oplus \RR$. If $P\in\calP'$ then
$P\Lambda P^{t}=\Lambda$, i.e.\ the true indeterminacy in the
choice of $\Lambda$ is measured by the Weyl 
group $W = \calP/ \calP'$, which in this case is just the full 
symmetric group $S_3$. (Note that $|\calP|=24$, and $|\calP'|=4$).
Again in analogy 
with three-body scattering, we think of the walls $\exp(\fraka) \setminus 
(\exp(\fraka))^*$ as `collision planes', where two eigenvalues are 
equal, and which we identify with the subset $\frakw \subset \fraka$. 
There are three collision planes, and these divide $\fraka$ into six chambers.

We next examine the structure of $M$ near the Weyl chamber walls.
Write the eigenvalues of $B \in M$, i.e.\ the diagonal entries
of $\Lambda$ in the decomposition for $B$ above, as
$\lambda_1$, $\lambda_2$ and $\lambda_3$. Suppose that $B$ 
lies in a small neighbourhood $\calU$ where 
\[
c< \lambda_1/\lambda_2< c^{-1},\quad \lambda_3>1/c,
\]
for some fixed $c \in (0,1)$. Recall also that $\lambda_3 =
1/\lambda_1 \lambda_2$. These inequalities imply that $\lambda_1=
(\lambda_1/\lambda_2)^{1/2}\lambda_3^{-1/2}<1$ and $\lambda_2= 
(\lambda_2/\lambda_1)^{1/2}\lambda_3^{-1/2}<1$, and $\lambda_3 > 1$
in $\calU$. Hence there is a well-defined decomposition 
$\RR^3 = E_{12} \oplus E_3$ for any $B \in \calU$, 
where $E_{12}$ is the sum of the first two eigenspaces and $E_3$ 
is the eigenspace corresponding to $\lambda_3$, regardless of whether 
or not $\lambda_1$ and $\lambda_2$ coincide. We could write equivalently
$B=O C O^t$, where $C$ is block-diagonal, preserving the splitting
of $\RR^3$. The ambiguity in this factorization is 
that $C$ can be conjugated by an element of $\On(2)$ (acting in the 
upper left block), and $\On(2)$ can be included in the top left
corner of $\SO(3)$ (the bottom right entry being set equal to $\pm 1$
appropriately). Let $C'$ denote the upper-left block of $C$; the
bottom right entry of $C$ is just $\lambda_3$, and so $\lambda_3 \det C' 
= 1$. In other words, $C' = \lambda_3^{-1/2}C''$ where $C''$ is
positive definite and symmetric with determinant $1$, hence
represents an element of $\SL(2)/\SO(2) \equiv \HH^2$. Hence for
an appropriate neighbourhood $\calV$ of $[I]$ in $\SL(2)/\SO(2)$, 
the neighbourhood $\calU$ is identified with $(\calV \times \SO(3))/\On(2) 
\times (1/c,\infty)$. Note that although the action of $\On(2)$ on $\calV$
has a fixed point (namely $[\Id]$), its action on $\SO(3)$,
and hence on the product, is free. The neighbourhood $\calV$
can be chosen larger when $\lambda_3$ is larger, and the limiting
`cross-section' $\lambda_3 = C$ has the form 
$(\HH^2 \times \SO(3))/\mbox{O}(2)$. This space is a fibre bundle
over $\SO(3)/O(2)$ ($= \RR P^2$) with fibre $\HH^2$. 
Notice that the Weyl chamber wall corresponds to the origin (i.e.\ 
the point fixed by the $\SO(2)$ action) in $\HH^2$. 

\subsection{Coordinates and metric}
We now discuss some useful coordinate systems on $M \setminus \calC$ using 
local coordinates on $\exp(\fraka)$ and $\SO(3)$. It is sufficient to
work near a diagonal element $B \in M$, and we write the diagonal
elements of elements in $\exp(\fraka)$ (near $B$) as $\lambda_1, 
\lambda_2, \lambda_3$ (in this order). Of course, these functions
are not independent since their product equals $1$. 
More general points in $M$ close to $B$ 
have $\SO(3)$ factor in their polar decomposition close to $\Id$, and 
so we augment these diagonal entries with the matrix entries 
$c_{12}, c_{13}$ and $c_{23}$ which parametrize the space of skew-symmetric 
matrices $\sol(3) = T_{\Id} \SO(3)$. Then a straightforward calculation
shows that at $B$, 
\begin{equation*}\begin{split}
g &= 6 \left( \frac{d\lambda_1^2}{\lambda_1^2} 
+ \frac{d\lambda_2^2}{\lambda_2^2} +\frac{d\lambda_3^2}{\lambda_3^2}\right)
\\
+ 3\left(\frac{\lambda_1}{\lambda_2} - \frac{\lambda_2}{\lambda_1}\right)^2
dc_{12}^2 
& + 3\left(\frac{\lambda_1}{\lambda_3} - \frac{\lambda_3}{\lambda_1}\right)^2
dc_{13}^2 
+ 3\left(\frac{\lambda_2}{\lambda_3} - \frac{\lambda_3}{\lambda_2}\right)^2
dc_{23}^2 
\end{split}\end{equation*}

Of course, a true coordinate system is obtained by fixing two
independent functions of the $\lambda_j$ as coordinates on $\exp(\fraka)$.
For reasons that will become clear later, an advantageous choice
is $\mu=\lambda_1/\lambda_2$ and $s = \lambda_3^{-3/2}$. We have
\[
\lambda_1=\mu^{1/2}\lambda_3^{-1/2} = \mu^{1/2}s^{1/3},\ 
\lambda_2=\mu^{-1/2}\lambda_3^{-1/2} = \mu^{1/2}s^{1/3},
\]
and hence 
\[
\lambda_2/\lambda_3=\mu^{-1/2}s, 
\ \lambda_1/\lambda_3=\mu^{1/2}s. 
\]
This gives
\begin{equation}
\begin{split}
g = & 3(d\mu/\mu)^2+ 4(ds/s)^2 + 
3(\mu-\mu^{-1})^2\,dc_{12}^2 
\\
+3s^{-2}(\mu^{1/2}s^2 & -\mu^{-1/2})^2\,dc_{13}^2
+3s^{-2}(\mu^{-1/2}s^2-\mu^{1/2})^2\,dc_{23}^2.
\end{split} \label{eq:met1}
\end{equation}

Notice that while this expression for the metric has 
coordinate singularities at the the Weyl chamber wall 
$\mu = 1$, $g$ is necessarily smooth across $\frakw$.
In fact, this is just a polar coordinate singularity.
To see this, write $\mu = e^r$, so $\mu=1$ corresponds to $r=0$ and
\begin{equation}
\begin{split}
& \qquad \qquad \quad g  = 3\left(dr^2 + \sinh^2 r\,d{\bar c}_{12}^2\right) 
\\ 
+ &\frac{4}{s^2}\left( ds^2 + \frac{3}{4}
(e^{r/2}s^2 -e^{-r/2})^2\,dc_{13}^2
+\frac{3}{4}(e^{-r/2}s^2-e^{r/2})^2\,dc_{23}^2\right),
\end{split} \label{eq:met1p}
\end{equation}
where $\bar{c}_{12} = 2 c_{12}$. It is clear that this is 
smooth across $r=0$ if $\theta = \bar{c}_{12} \in [0,2\pi]$,
and in fact the first part of this metric is just a multiple of 
the metric on hyperbolic space. 

We also wish to examine the metric $g$ in the region where 
\[
\mu \equiv \lambda_1/\lambda_2 < c, \qquad
\nu \equiv \lambda_2/\lambda_3 < c,
\]
for $c < 1$.  We have $\lambda_1 = \mu^{2/3}\nu^{1/3}$, $\lambda_2 = 
\mu^{-1/2}\nu^{1/3}$ and $\lambda_3 = \mu^{-1/3}\nu^{-2/3}$, and so 
\begin{equation}
\begin{split}
&g = 4( (d\mu/\mu)^2 + (d\mu/\mu)(d\nu/\nu) + (d\nu/\nu)^2)
\\
+ 3(\mu - &\mu^{-1})^2 dc_{12}^2 + 3(\mu\nu - \mu^{-1}{\nu}^{-1})^2
dc_{13}^2 + 3(\nu-\nu^{-1})^2 dc_{23}^2.
\end{split}
\label{eq:met2}
\end{equation}

The Laplacian in these coordinates is of the form
\begin{equation}
\begin{split}
\Delta_g = &\ \frac13\left((\mu D_\mu)^2+(\nu D_\nu)^2-(\mu D_\mu)(\nu D_\nu)
+i(\mu D_\mu)+i(\nu D_\nu)\right.\\
&\left.+(\mu D_{c_{12}})^2+(\nu D_{c_{23}})^2 
+(\mu\nu D_{c_{13}})^2\right) + E,
\end{split}
\label{eq:Lap-corner}
\end{equation}
where $E$ is the collection of all terms which are higher
order when $\mu$ and $\nu$ are small, i.e.\ it is a 
sum of up to two products of the vector fields $\mu D_\mu$,
$\nu D_\nu$, $\mu D_{c_{12}}$, $\nu D_{c_{23}}$ and $\mu \nu
D_{c_{13}}$, with at least one extra factor of $\mu$ or $\nu$.
We note that the `radial part' of this operator, which in
these coordinates simply corresponds to the operator 
acting on functions which are independent of the $c_{ij}$, is
\begin{equation}
\Delta_{\mathrm{rad}} = \frac13\left((\mu D_\mu)^2+(\nu D_\nu)^2-(\mu D_\mu)(\nu D_\nu)
+i(\mu D_\mu)+i(\nu D_\nu) + E'\right),
\label{eq:radLap}
\end{equation}
where $E'$ is an error term as above. In addition, using the
coordinates $\mu$ and $s$ which are valid in the interior of 
$H_\sharp$,
\begin{equation}
\begin{split}
\Delta_{\mathrm{rad}} = &
\frac13\left( (\mu D_\mu)^2 - \left(\frac{\mu + \mu^{-1}}{\mu-\mu^{-1}}
- \frac{s^2(\mu-\mu^{-1})}{s^4 - s^2(\mu+\mu^{-1})+1}\right) i \mu D_\mu\right)\\
&+ \frac14 \left( (sD_s)^2 - \frac{2(s^4-1)}{s^4 - s^2(\mu+\mu^{-1})+1}
i sD_s\right).
\end{split}
\label{eq:radLap2}
\end{equation}

\subsection{The compactification $\olM$}
We now proceed to describe the first compactification $\olM$ of $M$.
Consider first the interior of a Weyl chamber away from 
the walls. Fix the chamber where $\lambda_1 < \lambda_2 < \lambda_3$, 
and suppose that $0 < \mu,\nu <c<1$. As we have seen, $\SO(3)$ is just 
carried along as a factor in this region, and $(\mu,\nu,c_{12},c_{13},
c_{23})$ is a coordinate chart. We compactify by adding the faces
$\mu=0$ and $\nu=0$. These intersect at the corner $\mu=\nu=0$ at 
infinity.

Now consider the region $c < \mu < c^{-1}$, $\lambda_3 > c^{-1}$. 
We showed in \S 2.1 that this can be identified with 
$(\calV \times \SO(3))/\On(2) \times [1/c,\infty)$, where
$\calV$ is a ball in $\HH^2$. We compactify this by adding
the face $s' = 0$ where $s' \equiv \lambda_3^{-1}$, or in other 
words, we (partially) compactify this neighbourhood as $(\calV \times 
\SO(3))/\On(2) \times [0,c)_{s'}$. 

No new corners are added in this second step, and so the
induced compactification of $\fraka$ is as a hexagon, cf.\ 
the left picture in Figure~\ref{fig:flat2}. 

It remains to make sure that these two partial compactifications are 
compatible in the region of overlap, so that they fit together to 
give $\olM$ the structure of a compact $\Cinf$ manifold with 
corners of codimension $2$. In other words, we must show that the 
transition map, say when $c< \mu <1$, is smooth. But this is clear, 
for the transition map is given by diagonalizing the $2$-by-$2$ 
block $C''$ (cf.\ the discussion at the end of \S 2.1) and changing 
coordinates on the flat, which is certainly smooth away from the 
boundary at infinity. Notice, however, that the boundary defining
functions for this face ($\nu$ in the first chart and $s'$ in the second) are 
not smoothly related. In fact, since $\nu\mu^{1/2} = \lambda_3^{-3/2}$
and we are supposing that $\mu>c$ and $\nu \searrow 0$, we must
replace $s'$ by $s=\lambda_3^{-3/2}$ as the $\Cinf$ defining function
near this face. 

We summarize what we have done. The manifold $\olM$ is a compact 
manifold with a corner of codimension $2$. It has two boundary hypersurfaces,
$H_\sharp$ and $H^\sharp$, which are the closures of the parts of the 
boundary where the ratio of the smaller two, respectively the larger two, 
eigenvalues is bounded, or more directly, as the parts of the boundary 
where the ratio of the larger two, respectively the smaller two, 
eigenvalues vanishes. In the region $\mu, \nu \leq c < 1$, where
$\mu =\lambda_1/\lambda_2$, $\nu=\lambda_2/\lambda_3$, then 
\[
H_\sharp = \{\nu=0\}, \qquad H^\sharp = \{\mu=0\}.
\] 
The projection $((\HH^2 \times \SO(3))/\On(2)) \mapsto 
\SO(3)/\On(2)$ at $\lambda_3^{-1}=0$ is a fibration of the interior 
of $H_\sharp$ with base $\SO(3)/\On(2)$ and fibre $\HH^2$ and 
this extends to a fibration $\phi_\sharp$ of $H_\sharp$ with fibres 
$\overline{\HH^2} = \overline{B^2}$; there is an analogous fibration 
$\phi^\sharp: H^\sharp \to \SO(3)/\On(2)$. 

We remark that the metric $g$ induces a metric on the fibers of each of the faces. 
For example, from (\ref{eq:met1p}), we find that $g$ induces a metric
on $H_\sharp$, where $s=0$, which restricts to give ($3$ times) the 
standard hyperbolic metric on the $\HH^2$ fibres. 

\subsection{The boundary fibration structure}
This differential topological structure on $\olM$ is an example of a 
boundary fibration structure in the sense of \cite{RBM:Kyoto}. 
This particular boundary fibration structure, which we christen
the edge-to-edge (or ee) structure, is described in general as follows:

\begin{Def}
Suppose that $\overline{X}$ is a compact manifold with corners
of codimension $2$. Then we say that $\overline{X}$ is equipped
with an edge-to-edge structure if 

\begin{enumerate}
\item Each boundary hypersurface $H$ is the total space of a fibration
$\phi_H:H\to B_H$ which is transversal to $\pa H$ and has 
fiber $F_H$, a manifold with boundary;

\item Near any corner $H_1\cap H_2$, there is a fibre-preserving 
diffeomorphism of $\overline{X}$ to $F_{H_1}\times F_{H_2}\times \calU$, 
where $\calU$ is an open neighbourhood in some Euclidean space. 
\end{enumerate}

We write this structure as $(\overline{X},\phi)$, where $\phi$ 
is the collection of all of the mappings $\phi_H$. 
\end{Def}

Associated to an ee structure on $\overline{X}$ is the Lie algebra 
of all smooth vector fields on $\overline{X}$ which are arbitrary in 
the interior and are tangent not only to the boundaries $H$ (and hence 
to the corners $H_i \cap H_j$), but also the the fibres of the 
fibration $F_{H}$ along each codimension one boundary face $H$.  
These can be easily expressed in our particular case of interest, using the 
coordinates $\mu,\nu,c_{ij}$. The vector fields which meet all of 
these restrictions are in the span over $\Cinf(\olM)$ of the
basic vector fields
\begin{equation}
\mu \pa_\mu, \nu\pa_\nu,\ \mu \pa_{c_{12}},\ \nu \pa_{c_{23}},
\mu\nu \pa_{c_{13}}.
\label{eq:Vss-vf}
\end{equation}

\begin{Def}
If $(\overline{X},\phi)$ is an edge-to-edge structure,
then the associated Lie algebra of $\Cinf$ vector fields 
which are tangent to the fibers of $\phi_H$ for all $H$ is 
denoted by $\Vss(\overline{X})$. The enveloping algebra of this 
Lie algebra consists of all finite sums of products of elements
of $\Vss(\overline{X})$, and hence consists of a particular
class of differential operators, which we denote 
$\Diffss^*(\overline{X})$. 
\end{Def}

The space $\Vss(\bar X)$ is the full set of $\Cinf$ sections of a vector 
bundle $\Tss \overline{X}$ over $\overline{X}$, called the 
ee-tangent bundle. Its dual, the ee-cotangent bundle, is denoted 
$\Tss^*\overline{X}$. 

\begin{prop}
Let $\olM$ be the compactification of $M=\SL(3)/\SO(3)$ described above, 
and let $g$ be the invariant metric. Then $g \in \calC^\infty(\olM;S^2 (
\Tss^*\olM))$, and furthermore, $\Delta_g \in\Diffss^2(\olM)$.
\end{prop}

\section{The edge-to-edge small calculus}
The general strategy we shall follow, which has proved successful in 
many analogous situations, is that a boundary fibration structure leads 
to the introduction of an adapted calculus of pseudodifferential operators, 
and these operators in turn are used to investigate the analytic
properties of the elliptic operators associated to that boundary
fibration structure. As explained in the introduction, our aim here
is more modest in that we shall give a minimal development, including
only those parts of the theory of ee-pseudodifferential operators 
which are needed to understand the resolvent of $\Delta_g$ on the
space $M$. This simplifies the presentation considerably since it
amounts to only a little bit more than solving the `model problems'
which arise for more general operators (or for the Laplacian on
spaces endowed with an ee-metric, rather than the specific invariant 
one). We shall take up the more general theory elsewhere.

In this section we construct the `small calculus' of ee-pseudodifferential 
operators. These are characterized in terms of the behaviour of
their Schwartz kernels as conormal distributions on a resolution
of the double space $\overline{X}\times\overline{X}$. This is very
closely related to the construction of the calculus $\Psi_{p0}(\overline{X})$
in \cite{Mazzeo-Vasy:Resolvents}.

\subsection{The $\ssl$ double space}
We first describe the resolution process for $\olM$, where we have
already chosen good coordinates; after that we briefly describe the
process in the general situation. Thus the immediate goal is to
find a suitable blow-up $\olM^2_{\ssl}$ of $\olM^2$, which we call the 
edge-to-edge double space, where suitable means that the lift of
the Laplacian, first to either factor of $\olM^2$ and then to this
blow-up, is transversely elliptic with respect to the lifted diagonal, 
uniformly to all boundary faces.  Since $\Delta$ is an elliptic combination
of vector fields in $\Vss(\olM)$, this last property, in turn, will be
guaranteed provided the lifts of these vector fields span the normal
bundle of the lifted diagonal, uniformly on the closure of this
submanifold. (This is in contrast to what happens when these vector 
fields are lifted only to $\olM^2$, for on this space they all vanish at 
the boundary of the diagonal.) 

It will be simplest to describe the construction of $\olM^2_\ssl$
using local coordinates. 
We use two sets of coordinates, $(\mu,\nu,c_{ij},\mu',\nu',c_{ij}')$
on $\olM^2$ and consider especially the region near the product of
the corner with itself. We recall the simpler situation in the
edge calculus \cite{Mazzeo:Edge}, which is the special case of the 
ee-calculus when there is only one boundary hypersurface, which is
the total space of a fibration, and no corners. The appropriate
submanifold to blow up there is the fibre diagonal of the boundary,
and accordingly, away from the corner, e.g.\ where $\nu,\nu' \geq c > 0$
but $\mu,\mu' \to 0$, we should blow up the set 
\[
\calF^\sharp = \{\mu = \mu' = 0, c_{12} = c_{12}', c_{13}=c_{13}'\}.
\]
The effect of this blowup would be to desingularize the lifts
of the vector fields $\mu\del_\mu$, $\mu\del_{c_{12}}$ and
$\mu\nu\del_{c_{13}}$; the other vector fields $\nu\del_\nu$ and
$\nu\del_{c_{23}}$ are tangent to the submanifold being blown up and
are not affected. Similarly, in the region where $\mu,\mu' \geq c > 0$
but $\nu,\nu' \to 0$, we should blow up 
\[
\calF_\sharp = \{\nu = \nu' = 0, c_{23} = c_{23}', c_{13}=c_{13}'\}
\]
to desingularize the lifts of $\nu\del_\nu$, $\nu\del_{c_{23}}$ and
$\mu\nu\del_{c_{13}}$. Unfortunately, these two submanifolds 
do not intersect transversely, and to remedy this we must first
blow up their intersection 
\[
\calF^\sharp \cap \calF_\sharp = 
\{\mu=\mu'=\nu=\nu'=0, c_{ij}=c_{ij}' \ \ \forall i,j\}.
\]
The picture is slightly more complicated because to resolve
the quadratic singularity of the vector field $\mu\nu\pa_{c_{13}}$
it is necessary to perform a parabolic blow-up of this final
submanifold; this blow-up corresponds to the parabolic dilation
\[
\begin{split}
& (\mu,\mu',\nu,\nu',c_{12}-c_{12}', c_{12}', c_{23}-c_{23}', c_{23}',
c_{13}-c_{13}', c_{13}') \\ \qquad & \longmapsto  
(\lambda \mu,\lambda \mu', \lambda \nu,\lambda\nu',\lambda 
(c_{12} - c_{12}'), c_{12}', \lambda (c_{23}- c_{23}'), c_{23}', 
\lambda^2 (c_{13}-c_{13}'), c_{13}'), 
\end{split}
\]
for $\lambda > 0$. 

Now let us recast this more generally. Let $(\overline{X},\phi)$ 
be an ee structure. If $H$ is a boundary hypersurface, consider its 
fiber diagonal 
\[
\diag_{H,\phi_H}=\diag_H=\{(p,p')\in H\times H:\ \phi_H(p)=\phi_H(p')\}.
\]
For each (nontrivial) corner $K_{12}=H_1\cap H_2$, the submanifolds 
$\diag_{H_1}$ and $\diag_{H_2}$ intersect normally in a subset of 
$K_{12}\times K_{12}$. Thus the first step is to blow up each of these 
submanifolds $\diag_{H_1}\cap\diag_{H_2}$ parabolically, with the direction
of quadratic scaling given by the subbundle $S_{12} = 
(T\diag_{H_1}+T\diag_{H_2})\cap T_{\diag}$, 
where $\diag \subset \overline{X}^2$ is the diagonal. This makes 
the lifts of $\diag_{H_j}$ disjoint, and so we may now blow these
up in either order. Symbolically, if we enumerate the boundary
faces of $\overline{X}$ as $H_i$, then the ee double space of $\overline{X}$
is by definition the space
\[
\overline{X}^2_\ssl = \left[\overline{X}^2; \ \bigcup_{i,j} K_{ij}, \ 
\bigcup_{i,j} S_{ij};\ \bigcup_i \diag_{H_i}\right].
\]

\medskip

\noindent{{\it Coordinates on the blow-up.}\ } Each of these blow-ups 
can be realized by the introduction of appropriate polar coordinates. 
However, such coordinates are usually quite messy, and in practice it is 
far more convenient to work with projective coordinates as follows, even 
though these have the defect that they are singular at one of the 
boundary faces. Thus for example, let us introduce a set of
projective coordinates which is nonsingular away from the 
lift of the face where $\mu'=0$ (i.e.\ the copy of $H^\sharp$ on
the second copy of $\olM$). Then we may use
\begin{eqnarray*}
& \mu',\ \mu/\mu',\ \nu/\mu',\ \nu'/\mu', c_{12}',c_{13}',c_{23}' \\
& (c_{12}-c_{12}')/\mu', (c_{23}-c_{23}')/\mu', 
(c_{13}-c_{13}')/(\mu')^2
\end{eqnarray*}
as coordinates near the lift of $\calF^\sharp \cap \calF_\sharp$. 
Notice that in these coordinates $\mu'=0$ defines the `front
face', i.e.\ the new boundary hypersurface created in this blowup.
We continue by blowing up the lift of $\calF_\sharp$; in
these coordinates, 
\[
\calF_\sharp = \{\nu/\mu'=0,\ \nu'/\mu'=0,\ (c_{13}-c'_{13})/(\mu')^2=0,
\ (c_{23}-c'_{23})/\mu'=0\},
\]
and we shall use $\nu'/\mu'$ as the scaling coordinate. This gives
coordinates 
\begin{eqnarray*}
& \mu',\ s_1 = (\nu/\mu')/(\nu'/\mu') = \nu/\nu',\ 
s_2= \mu/\mu',\  \nu'/\mu', c'_{12},\ c'_{13},\ c'_{23}\\
& C_{12}=(c_{12}-c'_{12})/\mu',\ C_{23}=(c_{23}-c'_{23})/\nu',
C_{13}=(c_{13}-c'_{13})/(\mu'\nu').
\end{eqnarray*}
Here there are two boundary defining functions, $\mu'$ and
$\nu'/\mu'$; the former continues to define the front face from
the first blow-up, while the latter defines the new side-front
face introduced in this second blow-up. All other coordinates
are either angular or lateral along each of these faces.
To obtain coordinates near the lift of $\calF^\sharp$ it is 
most convenient to interchange the roles of $\mu'$ and $\nu'$
in the preceding discussion. 

\medskip

\noindent{{\it Vector fields.}\ } The virtue of using these projective
coordinates is immediately evident when we compute the lifts of the
basic vector field (\ref{eq:Vss-vf}) to $\olM^2_\ssl$. Thus, lifting
each of the vector fields (\ref{eq:Vss-vf}) first to the left
factor of $\olM$ in $\olM^2$ and then to the blowup gives
\begin{eqnarray*}
&\mu\del_\mu \to s_2\pa_{s_2},\  \nu\del_\nu \to s_1\pa_{s_1}, \\
& \mu\del_{c_{12}} \to s_2\pa_{C_{12}},\ \nu\del_{c_{23}} \to
s_1\pa_{C_{23}},\ \mu\nu\del_{c_{13}} \to s_1 s_2 \pa_{C_{13}}.
\end{eqnarray*}
The lifted diagonal $\diag_\ssl$ in $\olM^2\ssl$ is described 
by the equations $\{s_1 = s_2 = 1,\ C_{ij} = 0\}$, and from this 
it is clear that these vector fields span the normal bundle
to $\diag_\ssl$, uniformly to the boundary and corners. 

\subsection{$\ssl$-pseudodifferential operators}
Let $(\overline{X},\phi)$ be an ee structure, and let $\overline{X}^2_\ssl$
be the associated $\ssl$ double space. We now define the space
of $\ssl$-pseudodifferential operators $\Psi_{\ssl}^*(X)$ to
consist of those pseudodifferential operators $A$ on $X$ whose
Schwartz kernel $\kappa_A$ has the following properties. $\kappa_A$
is a distribution on $X^2$, and we require that it lift to
$\overline{X}^2_\ssl$ to be polyhomogeneous with respect to
the lifted diagonal $\diag_{\ssl}$, with singularities smoothly
extendible across all boundary faces of  $\overline{X}^2_\ssl$ 
which meet $\diag_\ssl$, and which vanishes to all orders
at all boundary faces which do not meet $\diag_\ssl$. For later
reference, we write $\calH'$ and $\calH''$ for the union of 
boundary faces which do or do not, respectively, meet $\diag_\ssl$. 

The preceding computations concerning the $\Vss$ vector fields 
immediately give the 
\begin{prop} If $L$ is any $\ssl$ differential operator, then
the Schwartz kernel of $L$ lifts to a (differentiated) delta section 
in $\overline{X}^2_{\ssl}$ supported along $\diag_\ssl$.
\end{prop}

As in ordinary pseudodifferential theory, there is a symbol map
\[
\sigma_m: \Psiss^m(\overline{X})\longmapsto S^m_{\hom}(\Tss^*\overline{X});
\]
an operator $L\in\Psiss^m(\overline{X})$ is called $\ssl$-elliptic 
if $\sigma_m(L)(z,\zeta)$ is invertible for $\zeta \neq 0$. 
Using this, the standard elliptic parametrix construction can be
mimicked to give
\begin{prop}
If $L\in\Psiss^m(\overline{X})$ is elliptic then there exists a parametrix
$G\in \Psiss^{-m}(\overline{X})$ such that $LG-\Id,GL-\Id\in\Psiss^{-\infty}
(\overline{X})$.
\end{prop}
 
Let us specialize again to the symmetric space $M$. For any $\ev \in \Cx$,
$\Delta - \ev$ is elliptic in this sense, and so this proposition yields
a parametrix $G(\ev) \in \Psiss^m(\overline{X})$ which depends 
holomorphically on $\lambda$, so that both $(\Delta - \ev)G(\ev)$
and $G(\ev)(\Delta - \ev)$ are of order $-\infty$ in this small
calculus. It is convenient to modify this parametrix slightly.
In fact, if $p\in M$ and $K_p$ is the stabilizer subgroup of $\SL(3)$ 
fixing this point, then $\Delta-\ev$ is $K_p$-invariant.
It would be nice to have a $K_p$-invariant parametrix, and this is
easy to arrange: simply define $G_p(\ev)=\int_{K_p} \phi_O^* G(\ev)
(\phi_O^{-1})^*\,dg_{K_p}(O)$, where $dg_{K_p}$ denotes the normalized 
invariant measure, and $\phi_O$ right multiplication by the element $O \in K_p$.
\begin{prop}
For any $\ev\in\Cx$ and $p\in M$, there exists an operator 
$G_p(\ev)\in \Psiss^{-2}(\overline{X})$ which is $K_p$-invariant and
which satisfies 
\begin{equation*}
(\Delta-\ev)G_p(\ev)-\Id,\ \ G_p(\ev)(\Delta-\ev)-\Id\in\Psiss^{-\infty}
(\overline{X})
\end{equation*}
Furthermore, $G(\ev)$ depends holomorphically on $\lambda$.
\end{prop}

As already explained in the introduction, this parametrix
is not the final one for the simple reason that the error terms
it leaves are not compact on $L^2$; removing these is the
principal motivation for defining the more elaborate large calculus. 

We conclude this section with a description of the regularity 
properties of this parametrix. We fix an $\ssl$-metric $g$, i.e.\ 
one for which the generating vector fields of the structure are of
bounded length. Boundedness of operators of order $0$ on $L^2(X,dV_g)$
may be deduced using the usual combination of a symbol calculus
argument to reduce to showing boundedness of operators of order
$-\infty$ and then proving this case directly by Schur's inequality;
cf.\ \cite{Mazzeo:Edge} for an example of this argument. Next,
for any $m \in \RR$, define the Sobolev space
\begin{equation}
\Hss^m(\olX) = \{u:  Au \in L^2(X;dV_g)\ \forall A \in \Psiss^m(\olX)\}.
\label{eq:sspee}
\end{equation}
If $m \in {\mathbb N}$, an equivalent formulation is
\[
\Hss^m(\olX) = \{u: V_1 \ldots V_j u \in L^2(X,dV_g),\ \text{
whenever}\ V_i \in \calV_{\ssl}\ \text{and}\ j \leq m\}.
\]
From the existence of the parametrix in the last proposition, we have
\begin{prop}
For any $m \in \RR$ and $\ev \in \Cx$, the parametrix $G_p(\ev)$ satisfies
\begin{eqnarray*}
& G_p(\ev):\Hss^m(\olM)\to \Hss^{m+2}(\olM), \\
& G_p(\ev):\Hss^m(\olM)^{K_p}\to \Hss^{m+2}(\olM)^{K_p}.
\end{eqnarray*}
In addition, if $E_p(\ev)=G_p(\ev)(\Delta-\ev)-\Id$ and
$F_p(\ev)=(\Delta-\ev)G_p(\ev) -\Id$ are the error terms, then 
\begin{equation*}
E_p(\ev),\, F_p(\ev):\Hss^m(\olM)^{K_p}\to \Hss^{\infty}(\olM)^{K_p}
\end{equation*}
for all $m$.
\label{pr:scp}
\end{prop}

To understand what these statements mean, suppose that $(\Delta - \ev)u = 
f \in \Cinf_0(M)$. If neither $u$ nor $f$ are $K_p$-invariant, then
from $u = G_p(\ev)f + E_p(\ev)u$ and assuming that $u \in L^2(M,dV_g)$, 
we get
\[
(\mu\del_\mu)^j(\nu\del_\nu)^k(\mu\del_{c_{12}})^r
(\nu\del_{c_{23}})^s (\mu\nu\del_{c_{13}})^t u \in L^2(M,dV_g)
\ \ \text{for all}\ j,k,r,s,t \geq 0.
\]
Thus $u$ has full tangential regularity only in directions tangent
to the fibres of the boundary fibration. On the other hand, suppose
that both $u$ and $f$ are $K_p$-invariant (and $f$ is still 
$\Cinf_0$). Then the vector fields $\del_{c_{ij}}$ annihilate $u$ 
and we see that $u$ restricts to a conormal function on the compactified
flat $\overline{\exp(\fraka)}$. In other words, 
\[
(\mu\del_\mu)^j(\nu\del_\nu)^k u \in L^2(M,dV_g)\ \ \text{for all}\ j,k \geq 0.
\]

\section{Product models}
The key remaining issue is to find a correction for the parametrix
$G_p(\ev)$ which solves away as much of the error terms $E_p(\ev)$
and $F_p(\ev)$ (as in Proposition \ref{pr:scp}) as possible. In fact, 
it will suffice to find a new parametrix for which the corresponding
error terms are not only regularizing, but which map $H^m_\ssl(\olM)^{K_p}$
into spaces of functions which decay at some definite rate
at infinity in the flat. As already indicated in the introduction,
this is done by solving away the expansions of the Schwartz kernels of 
these operators at the boundary faces $H_\sharp$ and $H^\sharp$
of $\overline{\exp(\fraka)}$. However, the main difficulty is caused
by the fact that the Laplacian does not have a product decomposition
at the corner $H_\sharp \cap H^\sharp$. This problem appears 
challenging, but as explained in detail in \S 5, the lift of
$\Delta_g$ to the logarithmic blow-up of $\olM$ is (quite remarkably)
well-approximated by product type operators at each of the new corners
of this blown up space. This means that the main technical difficulties
involve the analysis of the product operators arising in this process,
and this is the subject of the current section.  

We shall be using the results and methods of \cite{Mazzeo-Vasy:Resolvents}, 
concerning the detailed analysis of the resolvent of the Laplacian
on a product space $X = M_1 \times M_2$, which we review below.
Although that paper focused particularly on the case where both 
factors $(M_j,g_j)$ are conformally compact (asymptotically hyperbolic),
our application here requires that we let the factors be $(\Real^+_s,ds^2/s^2)$ 
and $(\HH^2,h)$, respectively (after all, the first is just the 
one-dimensional hyperbolic space). We shall both extend and refine the results 
in this setting. 

\subsection{Geometry and compactification of the product}
One of the main conclusions of \cite{Mazzeo-Vasy:Resolvents} is
that for purposes of analyzing its resolvent, or more precisely its Green
function with given pole, the best compactification
of $X=M_1\times M_2$ is given by 
\begin{equation}
\Xt=[(\bM_1)_{\log}\times(\bM_2)_{\log};\pa\bM_1\times\pa\bM_2].
\end{equation}
Recall that this means that if $\rho_j$, $j=1,2$, are smooth boundary
defining functions for $M_j$, then we replace these by $-1/\log \rho_j$ 
and then perform the standard blow-up of the corner. The function
\[
x = \sqrt{(\log \rho_1)^{-2} + (\log \rho_2)^{-2}}
\]
is a total boundary defining function for $\Xt$, i.e.\ is smooth
and vanishes simply at all the boundary faces with respect to this
new smooth structure. 

Now let
\[
(M_1,g_1) = \left(\RR^+_s, 4\frac{ds^2}{s^2}\right),\qquad
\mbox{and}\qquad 
(M_2,g_2) = (\HH^2,3h),
\]
where $h$ is the standard (curvature $-1$) metric on hyperbolic space.
If $\delta_j$ is the Riemannian distance function on each of these 
spaces, then the distance between pairs of points $(s,z), (s',z') 
\in X$ is given by 
\begin{equation*}
\delta((s,z),(s',z'))=\sqrt{\delta_1(s,s')^2+\delta_2(z,z')^2}.
\end{equation*}
In particular, fixing the point $o=(1,q) \in X$, then 
\[
\delta_1(s,1) = 2\, |\log s|, \qquad \delta_2(z,q) = \sqrt{3}\, |\log \mu|,
\]
where $\mu$ is a suitable defining function on $\HH^2$ (we include
the factors $2$ and $\sqrt{3}$ to keep track of the scaling factors 
on the metrics). 
Neglecting for the moment the fact that these functions are only
smooth away from $s=1$ and $z=q$, we set $\rho_1 = e^{-\delta_1(s,1)}$, 
$\rho_2 = e^{-\delta_2(z,q)}$. Hence 
\begin{equation}
x = \frac{1}{\delta(z,o)}=[4(\log s)^2+3(\log \mu)^2]^{-1/2}
\end{equation}
is a total boundary defining function for $\Xt$ and
\begin{equation}\begin{split}
&x_1=\frac{x}{(4(\log s)^2+1)^{-1/2}}\sim 2x\, |\log s|
=\frac{\delta_1(1,s)}{\delta(o,z)},\\
&x_2=\frac{x}{(3(\log \mu)^2+1)^{-1/2}}\sim \sqrt{3}x\, |\log \mu|
=\frac{\delta_2(q,z_2)}{\delta(o,z)}
\end{split}\end{equation}
are defining functions for the two side faces. Note that in terms of the
``eigenvalue coordinates'' $(\lambda_1,\lambda_2,\lambda_3)$ on $M$,
\[
x^{-2}=6((\log\lambda_1)^2+(\log\lambda_2)^2+(\log\lambda_3)^2).
\]

Returning to address the fact that the $\delta_j$ and hence $\delta$ 
are not smooth everywhere, we replace them by smoothed versions,
$\tilde\delta_j, \tilde\delta \in \Cinf$, which are chosen so that 
$\tilde\delta\geq 1$, 
$\delta\leq\tilde\delta\leq\delta+2$, and $\tilde\delta=\delta$ if
$\delta\geq 3$, and similarly for $\tilde{\delta}_j$. 
Although these no longer satisfy the triangle inequality, their failure to do
so is bounded: 
\begin{equation*}
\tilde\delta(z,z')+\tilde\delta(z',z'')-\tilde\delta(z,z'')\geq
\delta(z,z')+\delta(z',z'')-(\delta(z',z'')+2)\geq -2,
\end{equation*}
which is all we require in later estimates, and which we continue to
call the triangle inequality. Combining it with the fact that 
$\moddist\geq 1$, we have
\begin{equation*}
\moddist(z,z'')\leq \moddist(z,z')+\moddist(z',z'')+2
\leq 2\moddist(z,z')\moddist(z',z'')+2\leq 4\moddist(z,z')\moddist(z',z''),
\end{equation*}
and so 
\begin{equation}\label{eq:factor-bds}
\frac{\moddist(z,z'')}{\moddist(z,z')\moddist(z',z'')}\leq 4.
\end{equation}
We now replace our previous defining functions by 
\begin{equation*}
x_1=\frac{\moddist_1(1,s)}{\moddist(o,z)},\ x_2=\frac{\moddist_2(q,z_2)}
{\moddist(o,z)},\ x=\moddist(o,z)^{-1},
\end{equation*}
which are smooth and globally defined.

\subsection{Resolvent asymptotics}
The analysis of $\Delta_g$ on $K_p$-invariant functions 
ultimately reduces near the face $H_\sharp$ to that for the operator
\begin{equation*}
L_\sharp=\frac{1}{4}(sD_s)^2+i\frac{1}{2}(sD_s)
+\frac 13\Delta_h,
\end{equation*}
which is self-adjoint on 
\begin{equation*}
L^2(\Real^+_s\times\HH^2, s^{-3}ds\,dV_h).
\end{equation*}
It is computationally simpler to use its conjugate
\begin{equation}\label{eq:L_0-def}
L_0=s^{-1}L_\sharp s=\frac{1}{4}(sD_s)^2+\frac{1}{3}\Delta_h+\frac{1}{4},
\end{equation}
which is self-adjoint on 
\begin{equation*}
L^2(\Real^+_s\times\HH^2, dV_{g'}) = L^2(\Real^+_s\times\HH^2,
s^{-1}ds\,dV_h).
\end{equation*}
Note that $L_0 = \Delta_{g'} + \frac14$, where 
\begin{equation}\label{eq:L_0-metric}
g' = 4\frac{ds^2}{s^2}+3g_{\HH^2}.
\end{equation} 
(This explains our choice of factors $(M_j,g_j)$ above.) 

We now quote results from \cite{Mazzeo-Vasy:Resolvents} concerning
the structure of the resolvent 
\[
R_0(\ev) = (L_0 - \ev)^{-1} = (\Delta_{g'} -(\ev - \frac14))^{-1}.
\]
The formul\ae\ below frequently involve the quantity $\sqrt{\ev - \ev_0}$,
which is natural since 
\[
\inf\mbox{spec}\, (\Delta_{g'}) = 0 + \frac13 \cdot \frac14 = \frac{1}{12}
\Longrightarrow \inf\mbox{spec}\, (L_0) \equiv \ev_0 = \frac{1}{12} +
\frac14 = \frac13.
\]
We always use the branch of the square root which has {\em negative} 
imaginary part on $\Cx\setminus [0,\infty)$. 
These formul\ae\ also involve the Poisson operator $P_j(\ev)$, 
or rather its adjoint $P_j^t(\ev)$, and the spherical function $S_j(\ev)$,
on each of the factors. Of course, for the one-dimensional 
factor $M_1 = \RR^+$, these objects are particularly simple
and quite explicit. 

For simplicity, we begin with the asymptotics of $R_0(\ev)f$, where $f 
\in \dCinf(\Xb)$: 
\begin{prop}[\cite{Mazzeo-Vasy:Resolvents}, Proposition~7.7]
\label{prop:prod-asymp}
Let $f\in\dCinf(\Xb)$ and $\ev\in \Cx\setminus[\ev_0,\infty)$. 
Then on $\Xt$, 
\begin{equation}\label{eq:prod-asymp}
R_0(\mu)f=\mu^{1/2}x_2 \, x^{1/2}\exp (-i\sqrt{\ev-\ev_0}/x)\,h,\ \ 
h\in\Cinf(\Xt).
\end{equation}
Moreover, setting
\begin{equation}\label{eq:crit-point-im}
\ev_1^0(s)=\frac{\ev-\ev_0}{1+s^{2}},
\end{equation}
then the restriction of $h$ to the boundary is given by 
\begin{equation}\begin{split}\label{eq:pr-symbol-88}
&a(\ev,s)\left(P_1^t(\ev_1^0(s))\otimes P_2^t(\ev-\frac{1}{4}-
\ev_1^0(s))\right)f\ \text{on the front face},\\
& a'(\ev)\left(S_1(0)\otimes P_2^t(\ev-\frac{1}{4})\right)f\ 
\text{on the lift of}\ \bM_1\times \pa\bM_2,\\
&a''(\ev)\left(P_1^t(\ev-\ev_0)\otimes S_2(\frac{1}{12})\right)f\ 
\text{on the lift of}\ \pa\bM_1\times \bM_2, 
\end{split}\end{equation}
where $a,a'$ and $a''$ are all nonvanishing. 
\end{prop}

\begin{rem}
These formul\ae\ are different from the ones in \cite{Mazzeo-Vasy:Resolvents}
in one point: the resolvent $((sD_s)^2-\ev)^{-1}$ has a singularity of the
form $\ev^{-1/2}$ at the threshold branch point $0$, instead of one of
the form $\ev^{1/2}$ which occurs when this factor is higher dimensional. 
This simply causes the order of the leading term of the asymptotic 
expansions at the side faces where $s$ is finite to change by $1$,
and this produces an extra factor $x_2$, but not $x_1$, in 
\eqref{eq:prod-asymp}. In fact, the asymptotics at the lift of 
$\bM_1\times\pa\bM_2$ is even simpler than what we have stated 
since the $\Real^+$-invariance of $(sD_s)^2$ can be used to show that 
$R(\ev)f$ is still polyhomogeneous even when this face is blown down.
However, this does not affect the way we apply our results later, so we do
not take advantage of this.  
\end{rem}

From here, \cite{Mazzeo-Vasy:Resolvents} goes on to deduce the
full structure of the Schwartz kernel of the resolvent in the
product setting. Just as for the (more complicated) situation
we are in the process of establishing, that the resolvent for 
$\Delta$ on $M$ decomposes as a sum of two terms, $G_p(\ev)+ 
\tilde{R}(\ev)$, this kernel is a sum of two terms, $R_0'(\ev) + 
R_0''(\ev)$. The first term is in the small product-$0$ calculus 
(a simple case of the $\ssl$ calculus) and contains the full diagonal 
singularity, while the second is smooth on the interior but has 
a somewhat more complicated structure at the boundary which 
is resolved by passing to a further blow-up, the resolvent double space 
$\Xt^2_{\res}$. 

Operators in the small product-$0$ calculus $\Psi_{p0}(\Xb)$ (for
$X = M_1 \times M_2$), are characterized by the fact that their 
Schwartz kernels lift to the product-$0$ double space 
\begin{equation*}
\Xb^2_{p0} \equiv (\bM_1)^2_0 \times(\bM_2)^2_0
\end{equation*}
to be conormal to the lifted diagonal, smoothly extendible across
the front faces $\ff(\bM_1)^2_0 \times (\bM_2)^2_0$ and
$(\bM_1)^2_0 \times \ff(\bM_2)^2_0$ and vanishing to
infinite order at all other boundary faces. Here $M_1 = \Real^+$
and $M_2 = \HH^2$, so that $\bM_1 = I$ is the radial compactification 
of the half-line as an interval and $\overline{\HH^2}$ is the ball
$B^2$. 

The resolvent double space $\Xt^2_{\res}$ is obtained from $\Xb^2_{p0}$ 
by judiciously blowing up a certain (minimal) collection of corners, 
so that (at least when $\ev$ is in the resolvent set) $R_0''(\ev)$ lifts
to be polyhomogeneous. 

Recall that we have reduced to studying the restriction of the resolvent 
on $M$ to the flat; because of this, we only need to understand certain 
parts of the structure of the resolvent $R_0(\ev)$ for $L_0$ here. 
We have relegated the more detailed discussion of the structure of 
$\Xt^2_{\res}$ to the end of this paper in an appendix. 
More specifically, it clearly suffices to focus on the action of 
$R_0(\ev)$ on $\SO(2)$-invariant functions (with respect to a fixed point 
$q$) on the $\HH^2$ factor. We regard these functions as 
depending on the variables $s \in \RR^+$ and $\mu \in (0,1)$, 
where $\mu$ is the boundary defining function on $\overline{\HH^2}$
used earlier. The apparent `boundary' $\mu=1$ here is artificial,
and corresponds to the point $q$ in polar coordinates, and 
we systematically ignore it (for example by only considering
functions which are supported away from this set). In fact, set
\[
\frakb^+ = \Real^+_{s}\times(0,1)_{\mu} \subset \frakb = 
\Real_{s}^+ \times \Real_{\mu}^+. 
\]
If $\phi\in\Cinf(\RR^+)$, $\phi = 0$ near $\mu=0$ and $\phi=1$ for 
$\mu \geq 1/2$, then the Schwartz kernel of 
\begin{equation*}
(1-\phi(\mu))R_0(\ev)(1-\phi(\mu))
\end{equation*}
lifts to the resolvent double-space $\widetilde{\frakb}^2_{\res}$ 
and is supported on the closure of $\frakb^+\times\frakb^+$. 
Note that the Schwartz kernels of $\phi(\mu)R_0(\ev)(1-\phi(\mu))$
and $(1-\phi)(\mu)R_0(\ev)\phi(\mu)$ also lift trivially since
they are supported away from the corner $\mu=\mu'=0$. 

The advantage gained by this reduction is that the geometry of both 
the product-$0$ double space $\overline(\frakb)^2_{p0}$ and the 
resolvent double space $\overline{\frakb}^2_{\res}$ are simpler than 
when the second factor has dimension bigger than one. Thus the $0$-double 
space of $I = \overline{\RR^+}$ is obtained from $I^2$ by blowing 
up the boundary of the diagonal, $\pa\diag$, in $I^2$,  
$I^2_0 = [I^2; \pa \diag]$. This has the effect of separating the left
and right boundary faces. (Strictly speaking, this is not quite
true because they still intersect at the off-diagonal corners of
this square, but since all kernels are supported away from these
points, we studiously ignore this small untruth.) When $\dim M > 1$,
the left and right faces of the $0$-double space $M^2_0$ are no
longer separated, and necessitates some extra blow-ups in the 
resolvent double space. In any case, if $\lf_j$, $\rf_j$, $j=1,2$, are 
the left and right faces on the two factors $I^2_0$ in 
$\overline{\frakb}^2_{p0}$, then let $\calS$ be the collection of 
codimension two corners 
$\lf_1\times\lf_2$, $\lf_1\times\rf_2$, $\rf_1\times\lf_2$,
$\rf_1\times\rf_2$; by construction (and again neglecting the 
off-diagonal corners), these do not intersect. 
Now replace the defining functions 
$\rho_{\lf_j}$ and $\rho_{\rf_j}$ at these faces in each factor with
$\calR_{\lf_j} = -1/\log \rho_{\lf_j}$, $\calR_{\rf_j} = -1/\log \rho_{\rf_j}$.
Notice that we are not changing the defining function at the front
faces $\ff_j$.  The product-$0$ space with these defining functions is 
now denoted $\overline{\frakb}^2_{p0,\log}$. We finally define
\[
\widetilde{\frakb}^2_{\res} = [\overline{\frakb}^2_{p0}; \calS].
\]

The following theorem follows directly from \cite{Mazzeo-Vasy:Resolvents}, 
cf.\ the Appendix below. 

\begin{thm}\label{thm:int-kernel} If $\ev\in\Cx\setminus\spec(L_0)$, then 
\begin{equation}\label{eq:R-p-pp}
(L_0-\ev)^{-1} = R_0(\ev) =R_0'(\ev)+R''_0(\ev),
\end{equation}
where $R'_0(\ev)\in\Psi^{-2}_{p0}(\Xb)$ and 
\begin{equation}
\label{eq:F-def} 
R''_0(\ev)=(\rho_{\lf_2}\rho_{\rf_2})^{1/2} \exp(-i\sqrt{\ev-\ev_0}\ 
\moddist(z,z'))\, F(\ev),\qquad  \ev_0=\frac{1}{3}. 
\end{equation}
The kernel $F(\ev)$ on the right here is of the form
\[
F(\ev)=\moddist(z,z')^{-3/2}\, \moddist_2(z_2,z'_2) F'(\ev)
\frac{ds'}{s'}\,dV_h'
\]
where the primes on the density factor mean that it is pulled back
from the second ($z'$) factor. The function $F'(\ev)$ is 
bounded and smooth on $X \times X$ and lies in $\calC^0(\Xt\times X)$. 
In fact, if $\phi(\mu)\in\Cinf_0(\RR^+)$, $\phi = 1$ for $\mu$ near $1$, 
then for any $\ep>0$, 
\begin{equation}
e^{-\ep\moddist(z', o)}\phi(\mu)\, F'(\ev)\in\Cinf(\Xt;L^\infty(X)).
\end{equation}
(Here $z$ and $z'$ lie in $\Xt$ and $X$, respectively.) 
Of course, $F'$ is considerably more regular than this, but its $L^2$ 
mapping properties are determined by what is stated here. 

The restriction of this operator to $\SO(2)$-invariant functions 
in the second factor can be obtained by averaging $R_0''(\ev)$
with respect to Haar measure on $\SO(2)$, and can then be 
regarded as living on $\frakb^2$ as above. Writing 
$F'_0(\ev)$ as the average of $F'(\ev)$, then for all $\ep>0$, 
\begin{equation}
e^{-\ep\moddist(z',o)} (1-\phi(\mu))F'_0(\ev)\in 
\Cinf\left(\widetilde{\frakb^+};
L^2_\bl(\frakb^+; \frac{ds'}{s'}\frac{d\mu'}{\mu'}\right).
\end{equation}
A similar formula obtains when the first and second factors
are interchanged. 
\end{thm}

\subsection{Boundedness of the resolvent on weighted spaces}
We now prove some refined mapping properties of the resolvent $R_0(\ev)$ 
which are required later. These results reflect the fact that when $\ev$ 
is away from the spectrum, the Schwartz kernel of $R_0(\ev)$ has 
`off-diagonal' exponential decay of order $-\kappa$, where 
\[
\kappa = - \im\sqrt{\ev-\ev_0} > 0, \qquad \ev_0 = \frac13.
\]
For simplicity we phrase our theorems in terms of the decay
or asymptotics of $R_0(\ev)f$ for various classes of functions
$f$, rather than in terms of the structure of the kernel itself.

There are three results in this direction. The first states that 
if $|\alpha| < \kappa$ and $f$ decays (or grows) like $e^{\alpha/x}$
then so does $R_0(f)$. Next, if $f$ decays like $e^{\alpha/x}$ with
$\alpha < -\kappa$, then $R_0(\ev)f$ decomposes as a sum of two terms,
one decaying at the same rate as $f$ and another which has an
expansion, but decays only like $e^{-\kappa/x}$. Finally, 
we show that if $f$ decays exponential (at a rate $\alpha \in (-\kappa,
\kappa)$) in some sector, then $R_0(f)$ decays even faster
in disjoint sectors, with rate depending on the angle between
the two.  

\begin{prop} Suppose $|\alpha|<-\kappa$. Then
\[
R_0(\ev): e^{\alpha/x}H^m_{\ssl}(\Xb) \longrightarrow e^{\alpha/x}
H^{m+2}_{\ssl}(\Xb)
\]
is bounded, where we are using the $L^2$ measure with respect
to the metric $g'$ and the $\ssl$ Sobolev spaces defined in 
(\ref{eq:sspee}) (applied to the simple $\ssl$ structure on $X = \RR^+ 
\times \HH^2$). Moreover, if $\alpha= -\kappa$, then  
\[
R_0(\ev): x^p e^{\alpha/x}H^m_{\ssl}(\Xb) 
\longrightarrow x^{-p} e^{\alpha/x}H^{m+2}_{\ssl}(\Xb)
\]
for any $p>2$. (The restriction $p>2$ is not optimal, but suffices for our
later use.)
\end{prop}

\begin{proof} We shall only prove boundedness between weighted
$L^2$ spaces; the boundedness between Sobolev spaces, and
the gain of $2$ in $\ssl$ regularity, is a simple consequence
of the fact that $L_0$ can be applied on either the left 
or right, cf.\ Lemma~\ref{lemma:b-lift}, and arbitrary powers
of it can be commuted through. We also consider only the
case where $-\kappa < \alpha \leq 0$, since one need only
reverse the roles of $z$ and $z'$ to handle the case $0<\alpha<\kappa$.

The conclusion of the theorem is equivalent to the boundedness
of the mapping
\begin{equation*}
e^{-\alpha/x}R_0(\ev)e^{\alpha/x}\in\bop(L^2_{p0}(\Xb)).
\end{equation*}
It can be verified by direct calculation that the lift of the function 
\[
A_\alpha \equiv e^{\alpha/x' - \alpha/x}
\]
to $\olX^2_{p0}$ is smooth up to the front faces. We can assume that
the small-calculus part $R_0'(\ev)$ of the resolvent has support 
not intersecting any of the boundaries except the front faces, and
so $A_\alpha$ is smooth and bounded on its support. 
Thus we can focus on the term $R_0''(\ev)$, or rather $A_\alpha R_0''(\ev)$.  

It is convenient to replace the measure $dV_{g'} = s^{-1}\mu^{-2}\,
ds \,d\mu \,dy$, where $y$ is the angular (tangential) variable in
$\HH^2$, by $s^{-1}\mu^{-1}\, ds \,d\mu\, dy$, so we set
\begin{equation*}
L^2_b(\Xb)=L^2(X;\frac{ds\,d\mu\, dy}{s\mu})=\mu^{-1/2}L^2(\Xb; dV_{g'}).
\end{equation*}
An operator $A$ is bounded on $L^2(\Xb;dV_{g'})$ if and only
if $\mu^{-1/2}A\mu^{1/2}$ is bounded on $L^2_b(\Xb)$. 
Thus we must prove that
\begin{equation}\label{eq:b-kernel-8}
e^{\alpha/x'-\alpha/x}(\mu'/ \mu)^{1/2}R(\ev) \in\bop(L^2_b(\Xb)).
\end{equation}

From Theorem~\ref{thm:int-kernel}, with $x(z)=\moddist(z,o)^{-1}$,
then this conjugated kernel has the form
\[
e^{-\alpha(\moddist(o,z)-\moddist(o,z'))-i\sqrt{\ev-\ev_0}\moddist(z,z')}
e^{-(\moddist_2(z_2,z'_2)+\moddist_2(q,z'_2)-\moddist_2(q,z_2))/2}F(\ev)
\]
\[
=e^{-\alpha(\moddist(o,z)-\moddist(o,z')-\moddist(z,z'))}
e^{-(\moddist_2(z_2,z'_2)+\moddist_2(q,z'_2)-\moddist_2(q,z_2))/2}
e^{(-i\sqrt{\ev-\ev_0}+\alpha)\moddist(z,z')}F(\ev), 
\]
where $F$ is as in \eqref{eq:F-def}. Using the triangle inequality
to bound the exponents in the first two terms on the right,
we can rewrite this as 
\begin{equation*}
e^{-\gamma\moddist(z,z')}G,\ \qquad G\in L^\infty
(X\times X),\ \gamma= \kappa + \alpha > 0. 
\end{equation*}
This is an element of $L^\infty(X_{z};L^1_b(X_{z'})) 
\cap L^\infty(X_{z'};L^1_b(X_{z}))$, and hence the
conclusion follows from Schur's lemma. 

If $\alpha= -\kappa$ we can argue similarly, except 
now the kernel is rewritten as
\begin{equation*}
\moddist(z,o)^{-p} \moddist(z',o)^{-p}G,\qquad G\in L^\infty(X^2). 
\end{equation*}
Since $p>2$, this is integrable as before. 
\end{proof}

The next result concerns the behaviour of $R_0(\ev)f$ when 
$f \in e^{\alpha/x}L^2$ for some $\alpha < -\kappa$. As expected,
this function is the sum of two terms, the first decaying at
the same rate as $f$ and the second having the decay of a 
homogeneous solution to $L_0 u = 0$. We note, however, that
there are some subtleties in describing the precise regularity
of this second term; these arise already for the resolvent on
$\HH^2$ \cite{Mazzeo:Edge} (hence a fortiori on $\RR^+ \times \HH^2$).
More specifically, suppose that $f\in \mu^\gamma H^\infty_0
(\overline{\HH^2})$, where $\gamma >|\im\sqrt{\ev- 1/4}|$, or
in other words, $(\mu\del_\mu)^j(\mu\del_y)^\ell f \in \mu^\gamma 
L^2(dV_h)$ for all $j,\ell \geq 0$. Setting $u=(\Delta_{\HH^2}-\ev)^{-1}f$,
then the basic (small calculus) regularity result states that $u$ also 
lies in $\mu^\gamma H^\infty_0(\overline{\HH^2})$. However,
$u$ has no greater tangential regularity than $f$ itself, i.e.\  
we do not expect that $\del_y^\ell u \in \mu^\gamma L^2$ for $\ell > 0$
unless the same is true for $f$ too.  Returning to $X = \RR^+
\times \overline{\HH^2}$, we are fortunately spared these considerations 
because our main interest is when $f$ is $\SO(2)$ invariant.

In the following, let $H^m_{p0}(\Xb)$ denote the space
which we formerly called $H^m_\ssl$. We shall also use the logarithmically 
blown up single space 
\[
\Xt = [I_{\log}\times (\overline{\HH^2})_{\log};
\pa I\times\pa \overline{\HH^2}].
\]

\begin{prop}\label{prop:asymp} 
If $f\in e^{\alpha/x}H^m_{p0}(\Xb)$ for some $\alpha < -\kappa$, 
then 
\begin{equation*}
R(\ev)f= \mu^{1/2}x^{1/2} x_2
e^{-i\sqrt{\ev-\ev_0}/x}h+e^{\alpha/x}H^{m+2}_{p0}(\Xb);
\end{equation*}
In general, $h$ is (at least) continuous on $\Xt$, but if $f$ is 
$\SO(2)$-invariant, then $h$ is smooth on  $\Xt$.
In particular, if $m=\infty$, and $f$ is $\SO(2)$-invariant then 
\begin{equation*}
R(\ev)f= \mu^{1/2}x^{1/2}x_2
e^{-i\sqrt{\ev-\ev_0}/x}h',\quad  h'\in\Cinf(\Xt).
\end{equation*}
\end{prop}

\begin{proof}
We must show that
\begin{equation}\label{eq:b-kernel-8p}
e^{i\sqrt{\ev-\ev_0}/x} \mu^{-1/2}x^{-1/2}x_2^{-1}
R_0(\ev)\mu^{1/2}e^{\alpha/x}\in\bop(L^2_b(\Xb),C^0(\Xt)),
\end{equation}
and as in the preceeding proposition, we may immediately
replace $R_0$ by $R_0''$, since the small calculus contribution
$R_0'$ causes no difficulties. The kernel of this operator
then takes the form
\[
K''= e^{-i\sqrt{\ev-\ev_0}(\moddist(z,z')-\moddist(z,o))-
\kappa\moddist(z',o)}\] \[
 \times (e^{(\alpha+\kappa)\moddist(z',o)}\moddist(z',o)^{5/2})
\cdot e^{-(\moddist_2(z_2,z'_2)+\moddist_2(q,z'_2)-\moddist_2(q,z_2))/2}
\cdot \tilde F(\ev),
\]
where
\[
\tilde F(\ev)=x^{-1/2}x_2^{-1}\moddist(z',o)^{-5/2}F(\ev)
=\moddist(z,o)^{3/2}\moddist_2(z_2,q)^{-1}\moddist(z',o)^{-5/2}F(\ev).
\]

Now, dropping the singular measure, $\tilde F(\ev)$ is bounded since 
on the one hand,
\begin{equation*}
F(\ev)=\moddist(z,z')^{-3/2}
\moddist_2(z_2,z_2')F'(\ev),\ F'(\ev)\in
L^\infty(X\times X),
\end{equation*}
but we also know that
\begin{equation*}\begin{split}
&\moddist(z,o)^{3/2}\moddist_2(z_2,q)^{-1}\moddist(z',o)^{-5/2}
\moddist(z,z')^{-3/2}\moddist_2(z_2,z_2')\\
&=\left(\moddist(z,o)
\moddist(z',o)^{-1}\moddist(z,z')^{-1}\right)^{3/2}
\left(\moddist_2(z_2,q)^{-1}\moddist_2(z_2,z_2')\moddist(z',o)^{-1}\right)
\end{split}\end{equation*}
is bounded using \eqref{eq:factor-bds}.  In fact, using the triangle
inequality we even get
\begin{equation}\label{eq:K''-bd}
e^{(-\alpha-\kappa)\moddist(z',o)}\moddist(z',o)^{-5/2}
K''\in L^\infty(\tilde X\times \bar X).
\end{equation}
Each factor here is continuous on $\Xt\times X$ (note that we are
making no claims about continuity at infinity in the second factor). 
Indeed, for the main term $\tilde F$ this follows from the asymptotics
of $F(\ev)$, while for the other factors this is true when $z'$
lies in a compact set in $X$ since the triangle inequality gives
that both $\moddist_2(z_2,z_2')- \moddist_2(q,z_2)$ and 
$\moddist(z,z')-\moddist(z,o)$ are bounded continuous functions then. 
So we actually have
\begin{equation}\label{eq:K''-cont-o}
e^{-(\alpha+\kappa)\moddist(z',o)}K''\in C^0(\tilde X\times X).
\end{equation}
Since in fact $e^{(\alpha+\kappa)\moddist(z',o)}\moddist(z',o)^{-5/2}$
is continuous on $\tilde X\times\Xb$ and vanishes at 
$\tilde X\times \pa\Xb$, we deduce from \eqref{eq:K''-bd} and 
\eqref{eq:K''-cont-o} that 
\begin{equation}\label{eq:K''-cont}
e^{(\alpha+\kappa)\moddist(z',o)/2}K''\in C^0(\tilde X\times \Xb).
\end{equation}
Finally, $e^{(\alpha+\kappa)\moddist(z',o)/2}\moddist(z',o)^{-5/2}
\in L^2_b(\Xb)$, and this gives 
\begin{equation*}
K''\in C^0(\Xt;L^2_b(\Xb)), \qquad
\mbox{and so} \qquad  K'': L^2_b(\Xb)\to C^0(\Xt)
\end{equation*}
is continuous. This proves the first part of the proposition.

To continue, note that we would obtain full asymptotic expansions 
in general if it were true that $K''\in \Cinf(\Xt;L^2_b(\Xb))$. However, 
this is false on account of the behavior of $F(\ev)$ at the front face 
$I^2_0\times\ff(\overline{\HH^2})$. However, choosing a cutoff function 
$\phi$ as in Theorem~\ref{thm:int-kernel}, then both $\phi(\mu) K''$ and 
$K''\phi(\mu')$ do satisfy this smoothness criterion, and map into
functions with full expansions. It remains only to consider 
$(1-\phi(\mu))K''(1-\phi(\mu'))$.

Let us denote the corresponding kernel by $\tilde K''$, and 
work on the base space $\frakb^+$. The desired result now 
follows easily from Theorem~\ref{thm:int-kernel}, since that result
implies that if $P\in\Diff^m(\widetilde{\frakb^+})$ and $\ep > 0$, 
then $e^{-(\alpha+\kappa)\moddist(z',o)}e^{-\ep\moddist(z',o)}
P\tilde K''$ is continuous and bounded on $\widetilde{\frakb^+}
\times\frakb^+$. 

\end{proof}

The final result concerns off-diagonal decay of $R_0(\ev)$. 
To state this, choose $\SO(2)$-invariant cutoff functions 
$\phi$ and $\psi$ on the logarithmically blown up single space $\Xt$
with $\supp\phi\cap\supp\psi=\emptyset$. We can regard both as
defined on $\frakb_+$, or rather, its logarithmically blown up
single space $\widetilde{\frakb_+}$. As explained in the introduction,
the front face of $\widetilde{\frakb_+}$ can be identified with
(a large sector in) the Euclidean radial compactification 
$\overline{\frakb_r}$ of $\frakb_+$. This allows us to define
the angle $\theta$ between $\supp\phi$ and $\supp\psi$ on the
sphere at infinity. If $\supp\phi$ and $\supp\psi$ are conic outside 
a compact subset of $\fraka$, then this is simply the angle between
these cones. We have already shown that for $\ev\nin\spec(L_0)$, 
$\psi R(\ev)\phi$ is bounded on the spaces $e^{\alpha/x}L^2_{p0}(\Xb)$, 
$|\alpha| < \kappa$, but we now show that this cut off kernel
actually improves decay.  

Before proceeding, we note that this phenomenon is a very familiar one.
Consider the Laplacian $\Delta$ on $\Real^2_w$. The plane wave solutions 
of $(\Delta-\ev)u=0$ are those of the form $u(w)= e^{i\sqrt{\ev}\omega_0
\cdot w}$, $|\omega_0|=1$. (As usual, we use the branch of the square 
root function with negative imaginary part on $\CC \setminus \RR^+$.) 
Now write $w$ in polar coordinates as $r\omega'$.  Fix $\omega$ and 
let $\cos\theta = \omega \cdot \omega'$, so $\theta$ is the angle between 
the source of the plane wave and $w$. Then $|u| = 
e^{-\im\sqrt{\ev}\cos\theta r}$, or in other words, the exponential
rate of attenuation of $u$ is proportional to $\cos\theta$.
This effect is directly attributable to the structure of the 
resolvent $R_0(\ev) = (\Delta-\ev)^{-1}$ itself. Indeed, if 
$\chi\in\Cinf(\sphere^1)$, $\chi \equiv 1$ near $\omega_0$, then 
\begin{equation*}
u=\chi(\omega) e^{i\sqrt{\ev}\omega_0\cdot w}-R_0(\ev)
(\Delta-\ev)(\chi(\omega) e^{i\sqrt{\ev}\omega_0\cdot w}).
\end{equation*}
We refer to \cite{Froese-Herbst:Patterns} for an interesting discussion 
about a localization of this phenomenon (still in Euclidean space), 
concerning `decay profiles' of solutions of $(\Delta - \lambda)u = 0$ 
which are only defined in cones.  

When $\ev \in \Real^+$, one no longer obtains such decay, of course.
In its place is a propagation phenomenon at infinity, seen
already in \cite{RBMSpec}, which plays a very important
role in many-body scattering \cite{Vasy:Propagation-2}, 
\cite{Vasy:Bound-States}. The behaviour we are studying here,
which might be called ``dissipative propagation'', 
should be understood as a sort of analytic continuation
of these on-spectrum propagation results.

\begin{prop}[Dissipative propagation]
Choosing cut-off functions $\phi$, $\psi$ as above, let
$\theta\in(0,\pi/2)$ be less than the angle between their supports.
For $\ev\nin\spec(L_0)$ 
and $|\alpha|\leq\kappa=-\im\sqrt{\ev-\ev_0}$, write $\alpha=\kappa
\cos\theta_0$ for $\theta_0\in[0,\pi]$. Choose $\beta>\kappa\cos(
\theta+\theta_0)$ if $\theta+\theta_0\leq \pi$, otherwise choose 
$\beta>-\kappa$.  Then for any $m$,
\begin{equation*}
\psi R_0(\ev)\phi:
e^{\alpha/x}L^2_{p0}(\Xb)\to e^{\beta/x}H^m_{p0}(\Xb)
\end{equation*}
\end{prop}

\begin{proof}
Since we can assume that the supports of $\phi$ and $\psi$ do not
intersect, we can immediately discard the on-diagonal term
$R_0'(\ev)$, and furthermore it also suffices to just prove boundedness
into $e^{\beta/x}L^2_{p0}$. 

Let $\beta_0=\kappa\cos(\theta+\theta_0)$ if $\theta+\theta_0 \leq \pi$,
and otherwise let $\beta_0=-\kappa$. The key point is the
uniform boundedness of 
\begin{equation}\label{eq:exp-bded-16}
e^{\alpha\moddist(o,z')}e^{-\kappa\moddist(z,z')}
e^{-\beta_0\moddist(o,z)}
\end{equation}
when $z'\in\supp\phi$, $z\in\supp\psi$. It suffices to prove
the analogous boundedness when these modified distance functions
are replaced by the actual (nonsmooth) distance functions.
For $z=(s,z_2) \in \Real^+ \times \HH^2$, define
\[
w=(w_1,w_2)\in\Real^2,\ w_1=\log s=\delta_1(1,s),\ w_2=\delta_2(q,z_2),
\]
and define $w'$ analogously, corresponding to $z' = (s',z_2')$. Then
$\delta(o,z)= |w|$, $\delta(o,z')= |w'|$. The support conditions
on $\phi$ and $\psi$ mean that we restrict $w$ and $w'$ to lie
in cones $\Gamma$ and $\Gamma'$ in $\Real^2$  making an angle 
$\theta$ with one another.   
The triangle inequality gives $|\delta_2(z_2,o)-\delta_2(z_2',o)|\leq 
\delta_2(z_2,z_2')$, hence 
\begin{equation*}\begin{split}
\delta(z,z') & =\sqrt{\delta_1(s,s')^2+\delta_2(z_2,z_2')^2}
\geq\sqrt{\delta_1(s,s')^2+(\delta_2(z_2,o)-\delta_2(z_2',o))^2}\\
&=\sqrt{(w_1-w_1')^2+(w_2-w_2')^2}=|w-w'|.
\end{split}\end{equation*}
Hence the boundedness of \eqref{eq:exp-bded-16} follows from the estimate
\begin{equation}\label{eq:trig-est-8}
\alpha|w'|-\kappa|w-w'|-\beta_0|w|\leq 0, \qquad
\forall\ w \in \Gamma,\ w' \in \Gamma',
\end{equation}
or equivalently, dividing through by $\kappa>0$, 
\begin{equation}\begin{split}\label{eq:trig-est-16}
&\cos\theta_0|w'|\leq |w-w'|+ \cos(\theta+\theta_0)|w|, \qquad  \Mif
\theta+\theta_0\leq\pi,\\
&\cos\theta_0|w'|\leq |w-w'|-|w|, \qquad \qquad \qquad \ \  \, \Mif
\theta+\theta_0>\pi.
\end{split}\end{equation}
These certainly hold when one or the other vector vanishes. 
On the other hand, if neither is zero, then it suffices 
to consider the case where the angle between them is exactly
$\theta$, since this configuration minimizes the right hand side 
and keeps the left hand side fixed. 

Suppose first that $\theta+\theta_0>\pi$. Then $\pi\geq \theta_0>
\pi-\tilde\theta\geq 0$ implies $\cos\theta_0<\cos(\pi-\tilde\theta)=
-\cos\theta$, so we must only prove that
$|w|-\cos\tilde\theta|w'|\leq |w-w'|$. But this follows from
\[
|w|^2=(w-w')\cdot w+w'\cdot w\leq |w-w'||w|+|w'||w|\cos\tilde\theta
\]
once we divide through by $|w|$. 

If, on the other hand, $\theta+\theta_0\leq\pi$, then we let $v''$ 
be the unique unit vector such that $v'\cdot v''=\cos\theta_0$ and
$v\cdot v''=\cos(\theta+\theta_0)$. 
Since $w'\cdot v''=(w'-w)\cdot v''+w\cdot v''$
and $(w'-w)\cdot v''\leq |w-w'|$ as well, we deduce that $|w'|\cos\theta_0
\leq |w-w'|+|w|\cos(\theta+\theta_0)$, and this completes the proof.
\end{proof}

\begin{prop} With the same notation as above, if $f\in 
e^{\alpha/x}H^m_{p0}(\Xb)$ and $\theta+\theta_0>\pi$, then
\begin{equation*}
\psi R_0(\ev)\phi f= \mu^{1/2}x^{1/2} x_2
e^{-i\sqrt{\ev-\ev_0}/x}h+e^{\alpha/x}H^{m+2}_{p0}(\Xb),
\end{equation*}
where $h$ is continuous on $\Xt$. If $f$ is $\SO(2)$-invariant,
then $h\in\Cinf(\Xt)$. In particular, if $m=\infty$, and $f$ is 
$\SO(2)$-invariant, then
\begin{equation*}
\psi R_0(\ev)\phi f= \mu^{1/2}x^{1/2} x_2
e^{-i\sqrt{\ev-\ev_0}/x}h',\quad h'\in\Cinf(\Xt).
\end{equation*}
\end{prop}

\begin{proof}
We must now consider the kernel
\begin{equation*}\begin{split}
\tilde K''=&e^{(\alpha+\ep)\delta(o,z')}e^{-i(\sqrt{\ev-\ev_0}\delta(z,z')
-\delta(o,z))}\phi(z')\psi(z)\\
&\qquad\qquad\qquad\qquad
e^{-(\delta_2(z_2,z'_2)+\delta_2(q,z'_2)-\delta_2(q,z_2))/2}
\moddist(z',o)^{5/2}
\tilde F(\ev),
\end{split}\end{equation*}
where
\[
\tilde F(\ev)=x^{-1/2}x_2^{-1}\moddist(z',o)^{-5/2}e^{-\ep\delta(o,z')}F(\ev).
\]

Using \eqref{eq:exp-bded-16} with $\kappa$ in place of $\beta_0$, we deduce
that the first factor is bounded if $z\in\supp\psi$, $z'\in\supp\phi$,
and $\ep>0$ is sufficiently small.
Hence the arguments of
Proposition~\ref{prop:asymp} show that $\tilde K''\in C^0(\Xt;L^2_b(\Xb))$,
giving the first part of the conclusion. The rest follows as in
Proposition~\ref{prop:asymp}.
\end{proof}

Of course, $R_0(\ev)$ preserves $\SO(2)$-invariance, and if
$f$ is rotationally invariant, then $R_0(\ev)f$ also has
an expansion.

We conclude by translating this result back to the original operator
$L_\sharp$ and function space $\calH$ on which it is self-adjoint
by reintroducing the weight $s$:

\begin{cor}\label{cor:L_sharp}
The operator $(L_\sharp-\ev)^{-1}$ is bounded on $e^{\alpha/x}\calH$ for
$|\alpha| < \kappa$. Moreover, 
for $\phi$, $\psi$, $\theta$, $\alpha$, $\beta$ as above, 
\begin{equation*}
\psi (L_\sharp-\ev)^{-1}\phi:
 e^{\alpha/x}\calH\to e^{\beta/x}\calH.
\end{equation*}
If $\alpha < - \kappa$ and $f\in e^{\alpha/x}H^\infty_{p0}(\Xb)^{\SO(2)}$, then
$(L_\sharp-\ev)^{-1}f$ has a full asymptotic expansion on 
$\Xt$ of the form 
\begin{equation*}
(L_\sharp-\ev)^{-1}f= s\mu^{1/2}x^{1/2} x_2e^{-i\sqrt{\ev-\ev_0}/x}h',
\ h'\in\Cinf(\Xt).
\end{equation*}
\end{cor}

\section{Radial solutions and the final parametrix}
We now return to our main problem, and apply the results of the last 
section to finish the construction of a parametrix for $\Delta-\ev$
with compact remainder. As explained earlier, this requires
finding a correction term for the ($K_p$-invariant) 
small-calculus parametrix $G_p(\ev)$. Recall that the error terms 
$E_p(\ev)$ and $F_p(\ev)$ from 
Proposition~\ref{pr:scp} left after this first stage are $K_p$-invariant, 
hence can be regarded as acting on functions on $\fraka$. They are in 
fact residual elements of the $\ssl$ calculus on $\overline{\fraka}$. Thus, 
at the very least, we must solve equations of the form 
$(\Delta - \ev)u = f$, where $f$ is polyhomogeneous on $\overline{\fraka}$.
We now turn to this task. 

The idea is that near $H_\sharp$, $\Delta$ is well approximated by 
the product operator $L_\sharp$, and similarly it is well approximated
near $H^\sharp$ by $L^\sharp$. However, to make this precise 
we must pass to the logarithmically blown up space $\widetilde{\fraka}$
and localize there.  To motivate this, recall from \S 2.2 the 
coordinates $\mu$, $\nu$ near the corner on $\overline{\fraka}$
and the expression (\ref{eq:radLap}) for the restriction of $\Delta$ to 
$K_p$-invariant functions, i.e.\ to $\fraka$.   
Since $\Delta_{\mathrm{rad}}$ is not product type, even asymptotically, 
near the corner, we do not know a priori how to invert it. However,
working near $\nu = 0$, first note by (\ref{eq:radLap2}) that
\[
\Delta_{\mathrm{rad}} = \frac13\left( (\mu D_\mu)^2 
+ \frac{\mu + \mu^{-1}}{\mu - \mu^{-1}}i \mu D_\mu\right) + 
\frac14 (sD_s)^2 + \frac12 i s D_s + s^2 E,
\]
where $E = a \mu D_\mu + b s D_s + c$, $a,b,c \in \Cinf$. 
The first term in parentheses on the right
is just the radial part of the Laplacian on $\HH^2$ (with respect to the 
metric $3h$), and were we to have kept better track of the angular
derivatives, we would see the complete Laplacian on $\HH^2$ here.
Hence at least in the interior of $H_\sharp$, $\Delta - L_\sharp$ is small. 
Now consider the situation near the corner more carefully. The
coordinates $\mu$, $\nu$ are valid here, and we define the
change of variables $t = \mu$, $\bar{s} = \mu^{1/2}\nu$, so $s = \bar{s}$,
where $s$ is the notation used in the introduction. Then
\[
\mu D_\mu = t D_t + \frac12 \bar{s} D_{\bar{s}}, \quad \nu D_\nu = 
\bar{s}D_{\bar{s}},
\]
and inserting these into (\ref{eq:radLap2}) gives
\[
\Delta_{\mathrm{rad}} = \frac13 ( (t D_t)^2 + i t D_t ) 
+ \frac14 (\bar{s} D_{\bar{s}})^2 + \frac{i}2 \bar{s} D_{\bar{s}} + E,
\]
where $E$ is a different error term, consisting of sums of
of smooth multiples of the vector fields $t D_t$ and $\bar{s} D_{\bar{s}}$ 
with additional factors of $t$ or $\bar{s}$ and a polyhomogeneous function in 
these new coordinates vanishing at least to order one at the corner. 
Once again we see that the first term in parentheses is the radial part of the Laplacian 
and the main terms of this expression (i.e.\ omitting the error term) are
the same as for $L_\sharp$. 

This change of coordinates near the corner appears rather complicated,
but in fact it represents a smooth change of coordinates near the
entire closure of $\widetilde{H_\sharp}$ in $\widetilde{\fraka}$. 
To see this, let $\bar{\mu} = -1/\log \mu$, $\bar{\nu} = -1/\log \nu$,
$\tau = -1/\log t$, $\sigma = -1/\log \bar{s}$, from which
\[
\tau = \bar{\mu}, \qquad \sigma = \frac{\bar{\mu} \bar{\nu}}{\frac12 \bar{\mu}
+ \bar{\nu}}.
\]
The functions 
\[
r = \frac12 \bar{\mu} + \bar{\nu},\ \alpha = \frac{\frac12 \bar{\mu} - 
\bar{\nu}}
{\frac12 \bar{\mu}+ \bar{\nu}}, \qquad
r' = \tau + \sigma,\ \alpha' = \frac{\tau - \sigma}{\tau+\sigma}
\]
give two sets of coordinates on $\widetilde{\fraka}$ near this
face, and we have finally 
\[
r' = \frac12 r (1+\alpha)(3-\alpha), \qquad 
\alpha' = (1+\alpha)/(3-\alpha).
\]
This proves the claim that this coordinate change is smooth on
$\widetilde{\fraka}$ near the corner of $\widetilde{H_\sharp}$
since $(r,\alpha)$ and $(r',\alpha')$ are both smooth coordinate
systems there.

Altogether, noting in particular that $\bar{s} \leq \nu$ near the 
corner, we have just established that 
\[
\Delta-L_\sharp:\Hss^m(\bar M)^{K_p}\to \rho_\sharp\Hss^{m-1}(\bar M)
\]
near the closure of $H_\sharp$, and the proper venue for this
approximation is $\widetilde{\fraka}$. 

Now choose a smooth $K_p$-invariant partition of unity $\phi_\sharp+
\phi^\sharp+\phi_0=1$ on ${\widetilde M}$ such that $\supp\phi_\sharp$ 
and $\supp\phi^\sharp$ are disjoint from $\widetilde{H^\sharp} \cup 
\fraka^\sharp$ and $\widetilde{H_\sharp}\cup \fraka_\sharp$,
respectively, and $\supp\phi_0$ is compact in $M$. (We do not
introduce new notation for the lifts of $\fraka^\sharp$ and 
$\fraka_\sharp$ to $\widetilde{\fraka}$ because the walls of
the Weyl chambers are disjoint from where the blowup takes place.)
Also, choose additional cutoffs $\psi_\sharp$ and $\psi^\sharp$
which are identically $1$ on $\supp\phi_\sharp$ and $\supp\phi^\sharp$, 
respectively and with supports disjoint from $\widetilde{H^\sharp} 
\cup\fraka^\sharp$ and $\widetilde{H_\sharp} \cup \fraka_\sharp$. 

We now define
\begin{equation*}\begin{split}
\tilde R(\ev)=&G_p(\ev)-\psi^\sharp(L^\sharp-\ev)^{-1}\phi^\sharp
F(\ev)-\psi_\sharp
(L_\sharp-\ev)^{-1}\phi_\sharp F(\ev),
\end{split}\end{equation*}

Notation is abused slightly here because $L^\sharp$ is not
defined near $\fraka_\sharp$, but this is why the additional
cutoff $\psi^\sharp$ is included, and similarly for $L_\sharp$.

We compute that on $K_p$-invariant functions
\begin{equation}\begin{split}\label{eq:Rt-error}
& (\Delta-\ev)\tilde R(\ev)=\Id\, +\, (1-\phi_\sharp-\phi^\sharp)\, F(\ev) \\
& \qquad -(\Delta-L_\sharp)\, \psi_\sharp\, 
(L_\sharp-\ev)^{-1}\, \phi_\sharp \, F(\ev)
-[L_\sharp,\psi_\sharp]\, (L_\sharp-\ev)^{-1}\, \phi_\sharp \, F(\ev)\\
&\qquad -(\Delta-L^\sharp)\, \psi^\sharp\, (L^\sharp-\ev)^{-1}\,
\phi^\sharp \, F(\ev) -[L^\sharp,\psi^\sharp]\, (L^\sharp-\ev)^{-1}\,
\phi^\sharp \, F(\ev).
\end{split}\end{equation}

\begin{prop} The following maps are bounded
\begin{itemize}
\item[i)] For $|\alpha|<\kappa=-\im\sqrt{\ev-\ev_0}$,
\[
\tilde{R}(\ev): e^{\alpha/x}H^m_\ssl (M)^{K_p} 
\longrightarrow e^{\alpha/x}H^{m+2}_\ssl(M)^{K_p}
\] 
\item[ii)] For $\alpha< -\kappa$, 
\begin{equation*}
\begin{split}
\tilde R(\ev): e^{\alpha/x}\Hss^m(\bM)^{K_p} 
& \longrightarrow \\
e^{-i\sqrt{\ev-\ev_0}/x}\, x^{1/2}\, x_2\, \rho_\sharp\, \rho^\sharp\, 
\Cinf(\Mt) \ & +\ e^{\alpha/x}\Hss^{m+2}(\bM)^{K_p}. \end{split}
\end{equation*}
\item[iii)] In particular
\[ 
\tilde R(\ev): e^{\alpha/x}\Hss^\infty (\bM)^{K_p} 
\longrightarrow e^{-i\sqrt{\ev-\ev_0}/x}\, x^{1/2}\, x_2\, \rho_\sharp\,
\rho^\sharp\, \Cinf(\Mt).
\]
\end{itemize}
Now choose $\theta\in(0,\pi/2)$ less than the angles between
$\supp\phi_\sharp$ and $\supp d\psi_\sharp$ and between
$\supp\phi^\sharp$ and $\supp d\psi^\sharp$ on the sphere at infinity.
Then there exists a $\gamma > 0$ with  the following properties:
\begin{itemize}
\item[iv)] If $|\alpha|\leq \kappa$, so that $\alpha=\kappa
\cos\theta_0$ for some $\theta_0\in[0,\pi]$, and if in addition
$\beta>\kappa\cos(\theta+\theta_0)$ when $\theta+\theta_0\leq \pi$
and $\beta>-\kappa$ otherwise, then the error term $\tilde F(\ev)=
(\Delta-\ev)\tilde R(\ev)-\Id$ satisfies 
\[
\tilde F(\ev): e^{\alpha/x} H^m_\ssl (M)^{K_p}\longrightarrow  
e^{\max(\alpha-\gamma,\beta)/x}H^{m'}_\ssl(M)^{K_p}
\]
for any $m,m'$. 

\item[v)] If $\theta+\theta_0>\pi$ and $\alpha<-\kappa+\gamma$, then 
\[
\tilde F(\ev) :e^{\alpha/x}L^2(M)^{K_p} \longrightarrow 
e^{-i\sqrt{\ev-\ev_0}/x}\, x^{1/2}\, x_2\, \rho_\sharp\, \rho^\sharp\,
\Cinf(\Mt).
\]
\end{itemize}
\end{prop}

\begin{proof}
The mapping properties for $\tilde R(\ev)$ follow directly from
Corollary~\ref{cor:L_sharp}. As for the mapping properties for 
$\tilde F(\ev)$, note that the terms $[L_\sharp,\psi_\sharp]
(L_\sharp-\ev)^{-1}\phi_\sharp F(\ev)$, $[L^\sharp,\psi^\sharp]
(L^\sharp-\ev)^{-1}\phi^\sharp F(\ev)$ in \eqref{eq:Rt-error} are
bounded as stated since $[L_\sharp,\psi_\sharp]$ is supported on 
$d\psi_\sharp$, and so we may apply Corollary~\ref{cor:L_sharp}. 
On the other hand, to handle the remaining terms $(\Delta-L_\sharp)
\psi_\sharp 
(L_\sharp-\ev)^{-1}\phi_\sharp F(\ev)$, $(\Delta-L^\sharp)\psi^\sharp
(L^\sharp-\ev)^{-1}\phi^\sharp F(\ev)$, observe that 
$(L^\sharp-\ev)^{-1}\phi^\sharp F(\ev)$ is bounded on the appropriate
spaces, and then $(\Delta-L_\sharp)\psi_\sharp$ gives the additional
decay. This last statement follows from the fact that 
$\rho_\sharp\psi_\sharp < C\, e^{-\gamma/x}$ for some $\gamma>0$.
\end{proof}

The following theorem is now almost immediate.

\begin{thm}\label{thm:res-ee}
The resolvent $R(\ev)=(\Delta-\ev)^{-1}$ has the following mapping
properties on $K_p$-invariant functions. 
\begin{itemize}
\item[i)] For $|\alpha|<\kappa$, $m\in\Real$,
\begin{equation*}
R(\ev)\in\bop(e^{\alpha/x}\Hss^m(\bM)^{K_p},
e^{\alpha/x}\Hss^{m+2}(\bM)^{K_p});
\end{equation*}
\item[ii)] for $\alpha<-\kappa$, 
\begin{equation*}
R(\ev):e^{\alpha/x}\Hss^m(\bM)^{K_p} \longrightarrow 
e^{-i\sqrt{\ev-\ev_0}/x}\, x^{1/2}\, x_2\, \rho_\sharp\,\rho^\sharp\,
\Cinf(\Mt) +e^{\alpha/x}\Hss^{m+2}(\bM)^{K_p};
\end{equation*}
\item[iii)] in particular, 
\begin{equation}\label{eq:Green-asymp}
R(\ev):e^{\alpha/x}\Hss^\infty(\bM)^{K_p} \longrightarrow 
e^{-i\sqrt{\ev-\ev_0}/x}\,x^{1/2}\,x_2\,\rho_\sharp\,\rho^\sharp\, \Cinf(\Mt).
\end{equation}
\end{itemize}
\end{thm}

\begin{proof} 
We use the parametrix identity
\begin{equation}
R(\ev)=\tilde R(\ev)-\tilde R(\ev)\tilde F(\ev)
+\tilde E(\ev)R(\ev)\tilde F(\ev).
\label{eq:paramiden}
\end{equation}

Suppose first that $0\geq \alpha>-\kappa$. By the preceeding proposition, 
the first two terms are bounded $e^{\alpha/x}\Hss^m(M)^{K_p} \to 
e^{\alpha/x}\Hss^{m+2}(M)^{K_p}$. 

Next, write $\beta=\kappa\cos(\theta+\frac{\pi}{2})$ (corresponding to
$\cos\theta_0=0=0/\kappa$), and note that $\theta+\frac{\pi}{2}<\pi$.
We have $R(\ev), \tilde F(\ev) \in \bop(L^2(M)^{K_p})$; in
addition, if $\beta$ is given by the usual prescription but with 
$\alpha$ replaced by $0$, and if we set $\beta'=\max(-\gamma,\beta)<0$, then
$\tilde E(\ev):L^2(M)^{K_p}\to e^{\beta'/x}\Hss^{m'}(\bM)^{K_p}$. 
Hence for any $m'$, we deduce that 
\[
R(\ev): e^{\alpha/x}\Hss^m(\bM)^{K_p}\longrightarrow 
e^{\alpha'/x}\Hss^{m+2}(\bM)^{K_p}, \qquad \mbox{where}\quad
\alpha'=\max(\alpha,\beta').
\]

We now bootstrap to improve the image space. Thus suppose that
\begin{equation*}
\alpha_0=\inf\{\alpha'\in[\alpha,0]:
\ R(\ev)\in\bop(e^{\alpha/x}\Hss^m(\bM)^{K_p},\, 
e^{\alpha'/x}\Hss^{m+2}(\bM)^{K_p})\}.
\end{equation*}
We know that this set is non-empty, and that $\alpha_0<0$.
If $\alpha_0>\alpha$, then fix $\alpha'>\alpha_0$. Then by
definition of this inf we know that
\begin{equation*}
R(\ev)\in\bop(e^{\alpha/x}L^2(M)^{K_p}, e^{\alpha'/x}L^2(M)^{K_p});
\end{equation*}
Furthermore, both $\tilde R(\ev)$ and $\tilde F(\ev)$ are bounded on 
$e^{\alpha/x}L^2(M)^{K_p}$, while
\begin{equation*}
\tilde E(\ev)\in\bop(e^{\alpha'/x}L^2(M)^{K_p},e^{\alpha''/x}L^2(M)^{K_p}),
\qquad \mbox{where}\quad \alpha''=\min(\alpha'-\gamma,\beta').
\end{equation*}
Here $\beta'$ is calculated from $\alpha'$ as usual. Choosing
$\alpha'$ sufficiently close to $\alpha_0$, this gives
\begin{equation*}
R(\ev)\in\bop(e^{\alpha/x}L^2(M)^{K_p},
e^{\alpha''/x}L^2(M)^{K_p}),\qquad \mbox{where}\quad 
\alpha''<\alpha_0.
\end{equation*}
But this is a contradiction, and hence necessarily $\alpha_0=\alpha$. 
The same argument shows that $R(\ev)\in \bop(e^{\alpha/x}L^2(M)^{K_p})$.

Now suppose that $\alpha<-\kappa$. The first two terms
of (\ref{eq:paramiden}) map into 
\begin{equation}\label{eq:smooth-target-space-16}
e^{-i\sqrt{\ev-\ev_0}/x}x^{1/2}\,x_2\, \Cinf(\Mt)
+e^{\alpha/x}\Hss^{m+2}(\bM)^{K_p}.
\end{equation}
On the other hand,
\begin{equation*}
\tilde F(\ev)\in\bop(e^{\alpha/x}L^2(M)^{K_p},\, e^{\alpha'/x}L^2(M)^{K_p})
\end{equation*}
for any $\alpha'>-\kappa$, and at the same time, $\tilde R(\ev)
\in\bop(e^{\alpha'/x}L^2(M)^{K_p})$. But choosing
$\alpha'$ sufficiently close to $-\kappa$ ensures that
$\tilde E(\ev)$ maps this space to
\begin{equation*}
e^{-i\sqrt{\ev-\ev_0}/x}x^{1/2}x_2\Cinf(\Mt).
\end{equation*}
This finishes the proof.
\end{proof}

\section{Spherical functions}
The space of generalized eigenfunctions of $\Delta$ consists 
of the tempered distributions $u$ on $\olM$ such that $(\Delta-\ev)u=0$. 
If $p \in M$, then the $K_p$-invariant eigenfunctions are
called spherical functions. These were introduced by Harish-Chandra,
who analyzed them using a perturbation series. As noted in 
the introduction, the defect of this approach is that this series
only converges away from the Weyl chamber walls. However, 
it is straightforward to use the information obtained about the
resolvent to determine the behavior of these functions.

Suppose that $\xi\in\fraka_{\Cx}^*$ with $\im\xi\in
(\fraka^*)^+\equiv\fraka^+$ and $\xi\cdot\xi=\ev-\ev_0$ 
for some $\ev$ in the resolvent set $\Cx\setminus[\ev_0,\infty)$.
As always, set $\kappa = -\im \sqrt{\ev-\ev_0} > 0$. Then 
\begin{equation*}
u_0(z)=\rho_\sharp(z)\rho^\sharp(z)e^{-i\xi\cdot z}\in
e^{\kappa/x}x^r L^2(M)^{K_p},\ r<-1,
\end{equation*}
is an eigenfunction of $\Delta_{\mathrm{rad}}$. However, 
it is not smooth on $\Mt$ due to its behavior at the walls. 
Thus let $\psi$ be a cutoff function which vanishes near the walls, 
but is identically $1$ in an open cone around $\im\xi$.
(This is precisely where we are using that $\im\xi\in(\fraka^*)^+$.)
Then
\begin{equation*}
(\Delta-\ev)(\psi u_0)\in e^{\alpha/x} \Hss^\infty(\bM)^{K_p}
\end{equation*}
for some $\alpha<\kappa$, and we can thus use 
$R(\ev)$ to solve away this error term. We obtain the 
$K_p$-invariant eigenfunction
\begin{equation*}
U=\psi u_0-R(\ev)(\Delta-\ev)(\psi u_0),
\end{equation*}
which is precisely the spherical functions associated to the
parameter $\xi$. We have the preliminary result:
\begin{prop}
For $p\in M$, and $\xi\in\fraka_{\Cx}^*$ with $\im\xi\in(\fraka^*)^+$
and $\xi\cdot\xi=\ev-\ev_0$, there exists a $K_p$-invariant eigenfunction 
$U_\xi\in\Cinf(M)^{K_p}$ such that when $\alpha<\kappa$, 
\begin{equation}\label{eq:sph-fn-exist-8}
|U(z)-u_0(z)|\leq C\rho_\sharp(z)\rho^\sharp(z) e^{\alpha|z|},
\ z\in\overline{\fraka^+}.
\end{equation}
Moreover, $U_\xi$ is the unique $K_p$-invariant eigenfunction in 
$\dist(\bM)$ with 
\[
U(z)-u_0(z)\in\rho_\sharp(z)\rho^\sharp(z) e^{\alpha|z|}L^\infty(\fraka^+).
\]
\end{prop}

To obtain a finer description of the structure of $U = U_\xi$, we would
require a more explicit construction. This {\em can} be done 
were we to examine the parametrix construction more carefully.
However, we proceed as follows, relying on results from many-body 
scattering for motivation. In the present setting, this approach 
is essentially equivalent to the work of Trombi and Varadarajan
\cite{Trombi-Varadarajan:Spherical}, although it is more flexible
since it does not require the precise symmetric structure. 

In any case, by analogy with three-body scattering, we expect that 
$U$ decomposes as $u+v$, where $u$ is a sum of reflections of the initial 
plane wave $u_0$ and $v$ is an outgoing spherical wave, i.e.\ 
has the same asymptotics as the Green function \eqref{eq:Green-asymp}. 
Thus we posit that
\begin{equation}
u(z)=\rho_\sharp(z) \rho^\sharp(z)\sum_{\sigma\in S_3}c_\sigma(z)
e^{-i(\sigma\xi)\cdot z},
\label{eq:refl}
\end{equation}
for some choice of coefficients $c_\sigma$. We seek them in $\Cinf(\bM)$,
but note that this is much stronger than might be expected, since they
should presumably lie $\Cinf(\Mt)$ instead, which is indeed the
case in the analogous construction in three-body scattering.
We also demand that $c_1=1$ at the corner of $H_\sharp \cap H^\sharp$,
which ensures that each summand in (\ref{eq:refl}) agrees at the
corner.  The coefficients $c_\sigma$ are {\em not} well defined in 
the interior of the faces $H^\sharp$ and $H_\sharp$. Indeed, if
$\sigma \in S^3$ represents a reflection about one of the walls,
then for any $\xi'$, the functions $e^{-i(\sigma\xi')\cdot z}$ and 
$e^{-i\xi'\cdot z}$ have the same asymptotics at that wall.

Despite these concerns, we can determine the $c_\sigma$, which 
are constant along the main face, up to exponentially decreasing 
terms. The main requirement is that at each wall, the two reflected
terms associated to that wall must combine to yield the product of 
an exponential along the wall and a generalized eigenfunction of 
$\Delta_{\HH^2}$ in the normal direction. The $c_\sigma$ 
are all determined thereby from $c_1$ since reflections generate $S_3$. 
We now deduce that $(\Delta-\ev)u\in e^{\alpha/x}\Hss^\infty(\bM)^{K_p}$,
where we obtain the improvement in the weight to $\alpha$ by 
replacing $\Delta$ with $L_\sharp$ or $L^\sharp$ in the relevant regions;
in fact, we can take $\alpha=\kappa-1$. Hence if 
$\kappa$ is small (so $\ev$ is near to the spectrum), then we
can apply Theorem~\ref{thm:res-ee} to get
\begin{equation*}
\begin{split}
U & = \rho_\sharp(z) \rho^\sharp(z)\sum_{\sigma\in S_3}c_\sigma(z)
e^{-i(\sigma\xi)\cdot z} + v, \\
v \, \in & \, e^{-i\sqrt{\ev-\ev_0}/x}\, x^{1/2}\, x_2\, \rho_\sharp\, 
\rho^\sharp \, h, \qquad  h\in\Cinf(\Mt).
\end{split}
\end{equation*}

It is now apparent how the symmetric space structure plays a
very specific role. From Harish-Chandra's work, cf.\ 
\cite{Helgason:Groups}, it is known that along the main
face, the spherical wave is absent and $U$ simply equals the
main term $u$. The following theorem is therefore a
consequence of both the scattering theory (which gives the 
structure near the walls) and the algebraic structure (which 
eliminates the spherical wave). 

\begin{thm}
If $p\in M$ and $\xi$ is chosen $\fraka_{\Cx}^*$ with $\im\xi\in
(\fraka^*)^+$ and $\xi\cdot\xi=\ev-\ev_0$, $\ev\in\Cx\setminus[\ev_0,
\infty)$, then spherical function $U=U_\xi$ can be written as
\begin{equation}
U_\xi(z) =\rho_\sharp(z)\, \rho^\sharp(z)\, \sum_{\sigma\in S_3}
c_\sigma(z) e^{-i(\sigma\xi)\cdot z},\qquad
\mbox{where}\quad  c_\sigma\in\Cinf(\bM).
\end{equation}
\end{thm}

As we have already noted, this has been previously proved
by purely algebraic methods, cf.\ Trombi and Varadarajan, 
\cite{Trombi-Varadarajan:Spherical} and also Anker and Ji 
\cite{Anker-Ji:Heat}, Theorem~2.2.8.

Moreover, replacing $u_0$ by certain eigenfunctions of $L_\sharp$ or 
$L^\sharp$, we can also let $\xi$ approach the Weyl chamber 
walls; these special eigenfunctions of the model problems are
products of a radial solution of $(\Delta_{\HH^2}-\frac{1}{4})u_2=0$
(i.e.\ a spherical function at the threshold eigenvalue $1/4$)
with the exponential $u_1=s e^{i(-2\log s)\xi}$, 
$|\xi|^2=\ev-\ev_0$. The precise expression requires an 
additional factor involving the distance of $\xi$ from the collision 
plane.

\appendix

\section{Detailed description of the resolvent kernels}

For convenience, we review the definition of the resolvent space from 
\cite{Mazzeo-Vasy:Resolvents} which carries the Schwartz kernel
of the resolvent for the operator $L_0$ considered in \S 4. 
Recall that $X = \RR^+ \times \HH^2$, and we sometimes write
$M_1 = \RR^+$, $M_2 = \HH^2$, and also that $\overline{\Real_+} = I$.

The product $0$ double space $\Xb^2_{p0} = I^2_0 \times
(\overline{\HH^2})^2_0$ has six boundary hypersurfaces:
$$
\ff_1\times (\bM_2)_0^2,\ \lf_1\times (\bM_2)_0^2,\ \rf_1\times (\bM_2)_0^2,
\ (\bM_1)_0^2\times\ff_2,\ (\bM_1)_0^2\times\lf_2,\ (\bM_1)_0^2\times\rf_2;
$$
here $\ff_j$ is the front face, and $\lf_j$ and $\rf_j$ the left and right
faces of $(\bM_j)^2_0$.  We proceed from here by logarithmically
blowing up each of the two side faces $\lf_j$ and $\rf_j$ of
each factor, but {\it not} the front face. In other words, we just
introduce the new logarithmic defining functions 
\[
\calR_{\lf_j} = -1/\log \rho_{\lf_j},\qquad \calR_{\rf_j} = 
-1/\log \rho_{\rf_j}, \qquad j = 1,2,
\]
at these faces. The resulting space is denoted $(\bM_j)^2_{0,\log}$. 
We next blow up
\begin{equation}\label{eq:double-blow-up-1}
(\bM_1)^2_{0,\log}\times(\lf_2\cap\rf_2).
\end{equation}
In general we would also need to blow up $(\lf_1 \cap \rf_1) \times 
(\bM_2)^2_0$ now, but since $M_1$ is one dimensional,
$\lf_1\cap\rf_1 = \emptyset$.

The construction of the resolvent double space is completed by
blowing up the collection of submanifolds covering
\begin{equation}\begin{split}\label{eq:double-blow-up-2}
&(\lf_1\times (\bM_2)^2_{0,\log})\cap((\bM_1)^2_{0,\log}\times\lf_2),
\ (\lf_1\times (\bM_2)^2_{0,\log})\cap((\bM_1)^2_{0,\log}\times\rf_2),\\
&(\rf_1\times (\bM_2)^2_{0,\log})\cap((\bM_1)^2_{0,\log}\times\lf_2),
\ (\rf_1\times (\bM_2)^2_{0,\log})\cap((\bM_1)^2_{0,\log}\times\rf_2).
\end{split}\end{equation}
These are mutually disjoint after the blow-up in the previous step, 
and so these final blowups may be done in any order. 
Note that \eqref{eq:double-blow-up-1}-\eqref{eq:double-blow-up-2}
blow up {\em all five} codimension 2 corners coming from the intersections
of the four `side faces' of $\Xb^2_{p0}$:
$\lf_1\times (\bM_2)_0^2$, $\rf_1\times (\bM_2)_0^2$,
$(\bM_1)_0^2\times\ff_2$, $(\bM_1)_0^2\times\lf_2$, $(\bM_1)_0^2\times\rf_2$.

\begin{Def} The manifold with corners,
$\Xb^2_{\res}$, obtained by the series of blow-ups
of $\Xb^2_{p0}$ described above is called the resolvent compactification
of $X^2$.
\end{Def}

As proved in \cite{Mazzeo-Vasy:Resolvents}, the function
\begin{equation}\label{eq:S-def}
S=\frac{\delta_2(z_2,z_2')}{\delta_1(s,s')}
\end{equation}
and its inverse $S^{-1}$ are smooth on those portions of $\Xb^2_{\res}$ 
where they are bounded. In addition, 
\begin{equation*}\begin{split}
\calR=\delta(z,z')^{-1}&=\left(\delta_1(s,s')^2+
\delta_2(z_2,z_2')^2\right)^{-1/2}\\
&=\left((\calR_{\lf_1}^{-1}+R_{\rf_1}^{-1})^2
+(\calR_{\lf_2}^{-1}+\calR_{\rf_2}^{-1})^2\right)^{-1/2}
\end{split}\end{equation*}
is also smooth on $\Xb^2_{\res}$. 

\begin{lemma}\label{lemma:b-lift}
If $P\in\Diff_{p0}^m(\Xb)$ then $P$ lifts to elements
$P_L\in\Diff_b^m(\Xb^2_\res)$
and $P_R\in\Diff_b^m(\Xb^2_\res)$ from either the left or the right factors
of the projection $\Xb\times \Xb\to \Xb$.
\end{lemma}

Theorem~9.2 of \cite{Mazzeo-Vasy:Resolvents} describes the structure of
$(L_0-\ev)^{-1}$:

\begin{thm}\label{thm:res}
Suppose that $\ev\in\Cx\setminus\spec(L_0)$.
Then the Schwartz kernel of $R(\ev)=(L_0-\ev)^{-1}$
takes the following form:
\begin{equation}\begin{split}
&R(\ev)=R'(\ev)+R''(\ev),\ R'(\ev)\in\Psi^{-2}_{p0}(\Xb),\\
&R''(\ev)=(\rho_{\lf,2}\rho_{\rf,2})^{1/2}
\exp(-i\sqrt{\ev-\ev_0}/\calR) F(\ev),\ \ev_0=\frac{1}{3},
\end{split}\end{equation}
$F(\ev)$ is $\pi_R^*\Omega_{p0}\Xb$-valued polyhomogeneous
on $\Xb^2_{\res}$, order $0$ on the lift of the two front faces
$\ff_1\times(\bM_2)^2_{0,\log}$, $(\bM_1)^2_{0,\log}\times\ff_2$, order
$3/2$ on the lift of two of the `side faces' $\lf_1\times (\bM_2)_0^2$,
$\rf_1\times (\bM_2)_0^2$, as well as on the front face of the blow-up
\eqref{eq:double-blow-up-1},
and order $1/2$ on the other two `side faces'
$(\bM_1)_0^2\times\lf_2$, $(\bM_1)_0^2\times\rf_2$, as well as
on the front faces of the four blow-ups
\eqref{eq:double-blow-up-2}.
Moreover, the restriction of the leading term of $F(\ev)$ to the boundary,
i.e.\ its principal symbol, is
\begin{equation}\begin{split}\label{eq:pr-symbol-99}
&a(\ev,S)\left(P_1^t(\ev_1^0(S))\otimes P_2^t(\ev-\frac{1}{4}-\ev_1^0(S))\right)\ \text{on the front face of the blow-up of}\\
&\qquad\qquad\qquad\qquad\ (\lf_1\times (\bM_2)^2_{0,\log})
\cap((\bM_1)^2_{0,\log}\times\lf_2),\\
&
a'(\ev)\left(S_1(0)\otimes P_2^t(\ev-\frac{1}{4})\right)\ \text{on the lift of}\ (\bM_1)_0^2\times\lf_2,\\
&a''(\ev)\left(P_1^t(\ev-\ev_0)\otimes S_2(\frac{1}{12})\right)\ \text{on the lift of}\ \lf_1\times (\bM_2)_0^2,
\end{split}\end{equation}
with $\ev_1^0(S)$ given by
\begin{equation}\label{eq:crit-point-im-2}
\ev_1^0(S)=\frac{\ev-\ev_0}{1+S^{2}},
\end{equation}
$a,a',a''$ never zero, $S$ as in \eqref{eq:S-def}.
\end{thm}

In fact, we are mostly interested in $R(\ev)$ acting on $\SO(2)$-invariant
functions. Such functions can be regarded as living on
\begin{equation*}
\exp(-\frakb^+)=\Real^+_s\times
(0,1)_\mu,\ \frakb^+=\Real_{w_1}\times(0,\infty)_{w_2},
\ w_1=-\log s,\ w_2=-\log\mu;
\end{equation*}
the corresponding measure has a polar coordinate singularity at 
$\mu=1$. As indicated by the notation, it is sometimes convenient
to think of $\frakb^+$ as a subset of
\begin{equation*}
\frakb=\Real_{w_1}\times\Real_{w_2},
\end{equation*}
though we always restrict our attention to $w_2>0$.

So let $\phi\in\Cinf_c((0,\infty))$ be identically $1$
near $1$. Then the kernel of
\begin{equation*}
(1-\phi(\mu))R(\ev)(1-\phi(\mu))
\end{equation*}
can be described rather
simply on a compactification of $\frakb\times\frakb$, supported in
$\frakb^+\times\frakb^+$. More precisely, in complete analogy with
the construction above, we can consider the resolved double space
$\overline{\frakb^+}^2_{\res}$ of
\begin{equation*}
\overline{\frakb^+}=\overline{\Real^+_s}\times
[0,1)_\mu,
\end{equation*}
with the advantage that the
left and right faces from the same factor never intersect. Then it
is straightforward to describe the Schwartz kernel of
$(1-\phi(\mu))R(\ev)(1-\phi(\mu))$ as a distribution on
$\overline{\frakb^+}^2_{\res}$.

\begin{thm}\label{thm:res-flat}
Suppose that $\ev\in\Cx\setminus\spec(L_0)$ and
$\phi(\mu)\in\Cinf_c((0,\infty))$ is identically $1$
near $\mu=1$. Let
\begin{equation*}
\overline{\frakb^+}=\overline{\Real^+_s}\times
[0,1)_\mu.
\end{equation*}
Then the Schwartz kernel
of $(1-\phi(\mu))R(\ev)(1-\phi(\mu))$ acting on $\SO(2)$-invariant functions
takes the following form:
\begin{equation*}\begin{split}
&R'(\ev)+R''(\ev),\ R'(\ev)\in\Psi^{-2}_b(\overline{\frakb^+}),\\
&R''(\ev)=\mu^{1/2}(\mu')^{-1/2}\moddist(z,z')^{-3/2}\moddist_2(z_2,z_2')
\exp(-i\sqrt{\ev-\ev_0}/R)F_0(\ev),
\ \ev_0=\frac{1}{3},
\end{split}\end{equation*}
$F_0(\ev)$ is $\pi_R^*\Omega_b$-valued
smooth on $\overline{\frakb^+}^2_{\res}$.
\end{thm}

\begin{proof}
This follows by integrating out the $\SO(2)$ variables in the preceeding
theorem. The factor $(\mu')^{-1/2}$ appears as $(\mu')^{1/2}\pi_R^*\Omega_{p0}
=(\mu')^{-1/2}\pi_R^*\Omega_b$. However, it is instructive to see this 
directly, by constructing a parametrix for the restriction of $L_0$ 
to $\SO(2)$-invariant functions, 
which is an element of $\Diff^2_b(\overline{\frakb^+})$.
\end{proof}

\begin{proof}[Proof of Theorem~\ref{thm:int-kernel}]
\eqref{eq:R-p-pp}-\eqref{eq:F-def} follow immediately from
Theorem~\ref{thm:res-flat}.
\end{proof}

\bibliographystyle{plain}
\bibliography{sm}

\end{document}

%% file: macro.tex
\newcommand\Mand{\ \text{and}\ }
\newcommand\Mwith{\ \text{with}\ }
\newcommand\Mfor{\ \text{for}\ }
\newcommand\Mst{\ \text{such that}\ }
\newcommand\Mor{\ \text{or}\ }
\newcommand\Mif{\ \text{if}\ }
\newcommand\Miff{\ \text{iff}\ }
\newcommand\Mthen{\ \text{then}\ }
\newcommand\nin{\notin}
\newcommand\identity{\operatorname{id}}
\newcommand\Id{\operatorname{Id}}
\newcommand\Real{\mathbb{R}}
\newcommand\RR{\mathbb{R}}
\newcommand\CC{\mathbb{C}}
\newcommand\olM{\overline{M}}
\newcommand\olX{\overline{X}}
\newcommand\pos{\Real^+}
\newcommand\Rnp{\Real\setminus\{0\}}
\newcommand\nzero{\setminus\{0\}}
\newcommand\Cx{\mathbb{C}}
\newcommand\Cxp{\Cx^+}
\newcommand\Cxm{\Cx^-}
\newcommand\Nat{\mathbb{N}}
\newcommand\halfNat{{\frac{1}{2}}\mathbb{N}}
\newcommand\intgr{\mathbb{Z}}
\newcommand\im{\operatorname{Im}}
\newcommand\re{\operatorname{Re}}
\newcommand\sign{\operatorname{sign}}
\newcommand\codim{\operatorname{codim}}
\newcommand\End{\operatorname{End}}
\newcommand\Ker{\operatorname{Ker}}
\newcommand\Hom{\operatorname{Hom}}
\newcommand\tr{\operatorname{tr}}
\newcommand\Tr{\operatorname{Tr}}
\newcommand\ideal{{\mathcal I}}
\newcommand\Span{\operatorname{span}}
\newcommand\image{\operatorname{image}}
\newcommand\Range{\operatorname{Ran}}
\newcommand\Graph{\operatorname{graph}}
\newcommand\slim{\operatornamewithlimits{s-lim}}
\newcommand\gll{\mathfrak{gl}}
\newcommand\sll{\mathfrak{sl}}
\newcommand\sol{\mathfrak{so}}
\newcommand\GL{\operatorname{G\ell}}
\newcommand\SL{\operatorname{SL}}
\newcommand\SO{\operatorname{SO}}
\newcommand\On{\operatorname{O}}
\newcommand\pa{\partial}
\newcommand\del{\partial}
\newcommand\Rn{\Real^n}
\newcommand\Rm{\Real^m}
\newcommand\RN{\Real^N}
\newcommand\RtN{\Real^{2N}}
\newcommand\RM{\Real^M}
\newcommand\HH{\mathbb{H}}
\newcommand\sphere{\mathbb{S}}
\newcommand\Sn{\sphere^{n-1}}
\newcommand\Sm{\sphere^{m-1}}
\newcommand\Snp{\sphere^n_+}
\newcommand\Smp{\sphere^m_+}
\newcommand\SN{\sphere^{N-1}}
\newcommand\SNp{\sphere^N_+}
\newcommand\circlep{\sphere^1_+}
\newcommand\Phom{P_{h}}
\newcommand\Shom{S_{h}}
\newcommand\distance{\operatorname{dist}}
\newcommand\cl{\operatorname{cl}}
\newcommand\interior{\operatorname{int}}
\newcommand\Fa{\operatorname{Fa}}
\newcommand\ff{\operatorname{ff}}
\newcommand\mf{\operatorname{mf}}
\newcommand\cf{\operatorname{cf}}
\newcommand\scf{\operatorname{sf}}
\newcommand\lf{\operatorname{lf}}
\newcommand\rf{\operatorname{rf}}
\newcommand\indfam{{\mathcal K}}
\newcommand\ev{{\lambda}}
\newcommand\fraka{{\mathfrak a}}
\newcommand\frakb{{\mathfrak b}}
\newcommand\frakw{{\mathfrak w}}
\newcommand\calA{{\mathcal A}}
\newcommand\calB{{\mathcal B}}
\newcommand\calR{{\mathcal R}}
\newcommand\calS{{\mathcal S}}
\newcommand\calO{{\mathcal O}}
\newcommand\calJ{{\mathcal J}}
\newcommand\calM{{\mathcal M}}
\newcommand\calN{{\mathcal N}}
\newcommand\calX{{\mathcal X}}
\newcommand\calY{{\mathcal Y}}
\newcommand\calF{{\mathcal F}}
\newcommand\calG{{\mathcal G}}
\newcommand\calT{{\mathcal T}}
\newcommand\calC{{\mathcal C}}
\newcommand\calP{{\mathcal P}}
\newcommand\calU{{\mathcal U}}
\newcommand\calV{{\mathcal V}}
\newcommand\calCt{{\tilde {\mathcal C}}}
\newcommand\calCL{{\mathcal C}_{\text L}}
\newcommand\calCR{{\mathcal C}_{\text R}}
\newcommand\Cinf{{\mathcal C}^{\infty}}
\newcommand\dist{{\mathcal C}^{-\infty}}
\newcommand\dCinf{\dot{\mathcal C}^\infty}
\newcommand\ddist{\dot\dist}
\newcommand\Cj{{\mathcal C}^j}
\newcommand\Linf{L^{\infty}}
\newcommand\phg{{\text{phg}}}
\newcommand\bcon{{\mathcal A}}
\newcommand\bconc{{\mathcal A}_{\text{phg}}}
\newcommand\Sch{{\mathcal S}}
\newcommand\temp{\Sch^{\prime}}
\newcommand\Diff{\operatorname{Diff}}
\newcommand\Diffb{\operatorname{Diff}_{\text{b}}}
\newcommand\Diffc{\operatorname{Diff}_{\text{c}}}
\newcommand\Diffsc{\operatorname{Diff}_{\text{sc}}}
\newcommand\DiffI{\operatorname{Diff}_{\text{I}}}
\newcommand\DiffIq{\operatorname{Diff}_{\text{I},q}}
\newcommand\sing{\text{sing}}
\newcommand\reg{\text{reg}}
\newcommand\supp{\operatorname{supp}}
\newcommand\ssupp{\operatorname{sing\ supp}}
\newcommand\csupp{\operatorname{cone\ supp}}
\newcommand\esupp{\operatorname{ess\ supp}}
\newcommand\Fr{{\mathcal F}}
\newcommand\Frinv{\Fr^{-1}}
\newcommand\bop{{\mathcal B}}
\newcommand\spec{\operatorname{spec}}
\newcommand\pspec{\spec_{pp}}
\newcommand\cspec{\spec_{c}}
\newcommand\FIO{{\mathcal I}}
\newcommand\SP{\operatorname{RC}}
\newcommand\RC{\operatorname{RC}}
\newcommand\Symc{S_c}
\newcommand\Symca{S_c^{\alpha}}
\newcommand\Symczero{S_c^{0,...,0}}
\newcommand\sci{{}^{\text{sc}}}
\newcommand\sct{\sci T^*}
\newcommand\scdt{\sci \dot T^*}
\newcommand\dS{\dot S^*}
\newcommand\dT{\dot T^*}
\newcommand\dSreg{\dot\Sigma_{\text reg}}
\newcommand\scct{\sci\bar{T}^*}
\newcommand\Csc{C_{\text{sc}}}
\newcommand\SNpscd{(\SNp)^2_{\text{sc}}}
\newcommand\scdiag{\Delta_{\text{sc}}}
\newcommand\projscl{\pi^L_{\text{sc}}}
\newcommand\projscr{\pi^R_{\text{sc}}}
\newcommand\scHL{\sci H^{2,0}_{|\zeta|^2-\lambda^2}}
\newcommand\scHrg{\sci H^{2,0}_{\sqrt{g}}}
\newcommand\Hsc{H_{\text{sc}}}
\newcommand\WF{\operatorname{WF}}
\newcommand\WFp{\operatorname{WF^{\prime}}}
\newcommand\WFsc{\operatorname{WF}_{\text{sc}}}
\newcommand\WFscp{\operatorname{WF_{sc}^{\prime}}}
\newcommand\WFC{\operatorname{WF}_C}
\newcommand\WFCi{\operatorname{WF}_{C_i}}
\newcommand\elliptic{\operatorname{ell}}
\newcommand\Psop{\operatorname{\Psi}}
\newcommand\Psiscrs{\operatorname{\Psi_{sc}^{-2,\infty}}}
\newcommand\Psiscr{\operatorname{\Psi_{sc}^{-2,0}}}
\newcommand\Psiscrm{\operatorname{\Psi_{sc}^{0,2}}}
\newcommand\PsiscHam{\operatorname{\Psi_{sc}^{2,0}}}
\newcommand\Psisci{\operatorname{\Psi_{sc}^{*,*}}}
\newcommand\Psiscid{\operatorname{\Psi_{sc}^{0,0}}}
\newcommand\Psiscis{\operatorname{\Psi_{sc}^{0,\infty}}}
\newcommand\Psiscsi{\operatorname{\Psi_{sc}^{-\infty,0}}}
\newcommand\Psiscs{\operatorname{\Psi_{sc}^{-\infty,\infty}}}
\newcommand\Psiscalg{\operatorname{\Psi_{sc}^{\infty,-\infty}}}
\newcommand\nullHam{{\mathcal N}}
\newcommand\charD{\Sigma_{\Delta-\lambda^2}}
\newcommand\charLap{\Sigma_{\Delta-\lambda}}
\newcommand\Snl{\Sn_{\lambda}}
\newcommand\SNl{\SN_{\lambda}}
\newcommand\gammat{\tilde\gamma}
\newcommand\gammasc{\gamma}
\newcommand\Tau{\mathcal{T}}
\newcommand\taut{\tilde\tau}
\newcommand\taub{\bar\tau}
\newcommand\Nout{N^+_{\lambda}}
\newcommand\Nin{N^-_{\lambda}}
\newcommand\Nio{N^{\pm}_{\lambda}}
\newcommand\El{E_{\lambda}}
\newcommand\Elt{\tilde E_{\lambda}}
\newcommand\Eil{E^i_{\lambda}}
\newcommand\Ejl{E^j_{\lambda}}
\newcommand\Eajl{E^{\alpha_j}_{\lambda}}
\newcommand\Eilt{\tilde E^i_{\lambda}}
\newcommand\Np{N^+}
\newcommand\Nm{N^-}
\newcommand\Npm{N^{\pm}}
\newcommand\Fin{F^-(\lambda)}
\newcommand\Fini{F^-_i(\lambda)}
\newcommand\Fout{F^+(\lambda)}
\newcommand\Fouti{F^+_i(\lambda)}
\newcommand\Foutj{F^+_j(\lambda)}
\newcommand\Rout{R^+_{\lambda}}
\newcommand\Routl{R^+_{\lambda^2}}
\newcommand\Routsgnl{R^{\sign\lambda}_{\lambda^2}}
\newcommand\Rin{R^-_{\lambda}}
\newcommand\Rinl{R^-_{\lambda^2}}
\newcommand\Rinsgnl{R^{-\sign\lambda}_{\lambda^2}}
\newcommand\Rio{R^{\pm}_{\lambda}}
\newcommand\Riol{R^{\pm}_{\lambda^2}}
\newcommand\Roi{R^{\mp}_{\lambda}}
\newcommand\Roil{R^{\mp}_{\lambda^2}}
\newcommand\Riob{R^{\pm}}
\newcommand\Roib{R^{\mp}}
\newcommand\Tio{T^{\pm}}
\newcommand\Tiob{T^{\pm}_{\ff}}
\newcommand\Toi{T^{\mp}}
\newcommand\Toib{T^{\mp}_{\ff}}
\newcommand\TIiob{T_I^{\pm}}
\newcommand\Rinb{R^-}
\newcommand\Rinbsgnl{R^{-\sign\lambda}}
\newcommand\Tin{T^-}
\newcommand\Tinb{T^-_{\ff}}
\newcommand\TIinb{T^-_I}
\newcommand\Routb{R^+}
\newcommand\Routbsgnl{R^{\sign\lambda}}
\newcommand\Tout{T^+}
\newcommand\Toutb{T^+_{\ff}}
\newcommand\TIoutb{T^+_I}
\newcommand\Rlkf{(|\xib|^2-(\lambda-i0)^2)^{-1}}
\newcommand\Rlk{\rho_0(\lambda)}
\newcommand\Rmlk{\rho_0(-\lambda)}
\newcommand\Rpmlk{\rho_0(\pm\lambda)}
\newcommand\Rlka{\rho_1(\lambda)}
\newcommand\Rlkb{\rho_2(\lambda)}
\newcommand\Rilk{\rho_i(\lambda)}
\newcommand\reduced{\natural}
\newcommand\Rlf{R_0(\lambda)}
\newcommand\Rla{R_1(\lambda)}
\newcommand\Rlb{R_2(\lambda)}
\newcommand\Ril{R_i(\lambda)}
\newcommand\Rlj{R_j(\lambda)}
\newcommand\Rlft{R_0(\lambda)}
\newcommand\Rflambda{R_0^{\reduced}(\sigma)}
\newcommand\RV{R^{\reduced}_V}
\newcommand\Rfsigma{R_0^{\reduced}(\sigma)}
\newcommand\Rfsigmah{R_0^{\reduced}(\sigma^{1/2})}
\newcommand\Rfzero{R_0^{\reduced}(0)}
\newcommand\RlV{R^{\reduced}_V(\sigma)}
\newcommand\RlVi{R^{\reduced}_{V_i}(\sigma)}
\newcommand\RlVt{R_V(\lambda)}
\newcommand\RlVtL{{R}_V^L(\lambda)}
\newcommand\RlVtR{{R}_V^R(\lambda)}
\newcommand\RlVit{{R}_{V_i}(\lambda)}
\newcommand\RlVta{{R}_V^{(1)}(\lambda)}
\newcommand\RlVtk{{R}_V^{(k)}(\lambda)}
\newcommand\RlVatV{{R}_{V_{\alpha}}(\lambda)V_{\alpha}}
\newcommand\RlVatVa{{R}_{V_{\alpha_1}}(\lambda)V_{\alpha_1}}
\newcommand\RlVatVb{{R}_{V_{\alpha_2}}(\lambda)V_{\alpha_2}}
\newcommand\RlVatVk{{R}_{V_{\alpha_k}}(\lambda)V_{\alpha_k}}
\newcommand\RlVatVkk{{R}_{V_{\alpha_{k+1}}}(\lambda)V_{\alpha_{k+1}}}
\newcommand\RlVaptV{{R}_{V_{\alpha'}}(\lambda)V_{\alpha'}}
\newcommand\RlVapptV{{R}_{V_{\alpha''}}(\lambda)V_{\alpha''}}
\newcommand\RlVajtV{{R}_{V_{\alpha_j}}(\lambda)V_{\alpha_j}}
\newcommand\RlVaktV{{R}_{V_{\alpha_k}}(\lambda)V_{\alpha_k}}
\newcommand\RlVakktV{{R}_{V_{\alpha_{k+1}}}(\lambda)V_{\alpha_{k+1}}}
\newcommand\Tl{T(\lambda)}
\newcommand\Tlt{\tilde\Tl}
\newcommand\Tltp{\tilde T'(\lambda)}
\newcommand\Tltpp{\tilde T''(\lambda)}
\newcommand\Tli{T_i(\lambda)}
\newcommand\Tlit{\tilde\Tli}
\newcommand\Tlip{T_i'(\lambda)}
\newcommand\Tlipp{T_i''(\lambda)}
\newcommand\Tlj{T_j(\lambda)}
\newcommand\Tla{T_{\alpha}(\lambda)}
\newcommand\Tlaa{T_{\alpha_1}(\lambda)}
\newcommand\Tlab{T_{\alpha_2}(\lambda)}
\newcommand\Tlak{T_{\alpha_k}(\lambda)}
\newcommand\Tlakt{\tilde\Tlak}
\newcommand\Tlaj{T_{\alpha_j}(\lambda)}
\newcommand\Tlajj{T_{\alpha_{j+1}}(\lambda)}
\newcommand\Tlajp{T_{\alpha_j}'(\lambda)}
\newcommand\Tlajpt{\tilde\Tlajp}
\newcommand\Tlajt{\tilde\Tlaj}
\newcommand\Tlakk{T_{\alpha_{k+1}}(\lambda)}
\newcommand\Tlakkp{T_{\alpha_{k+1}}'(\lambda)}
\newcommand\Tlap{T_{\alpha'}(\lambda)}
\newcommand\Tlapt{\tilde\Tlap}
\newcommand\Tlapp{T_{\alpha''}(\lambda)}
\newcommand\Tkl{T^{(k)}(\lambda)}
\newcommand\Tcl{T^{\flat}(\lambda)}
\newcommand\Fl{F(\lambda)}
\newcommand\BlVt{\tilde B_V(\lambda)}
\newcommand\KBlVt{K_{\BlVt}}
\newcommand\BlVaat{B_{V_{\alpha_1}}(\lambda)}
\newcommand\BV{B_V}
\newcommand\Bone{B_1}
\newcommand\Btwo{B_2}
\newcommand\Bthree{B_3}
\newcommand\Banyj{B_j}
\newcommand\PlV{P_V(\lambda)}
\newcommand\PlVc{P_V^{\flat}(\lambda)}
\newcommand\Pl{P_0(\lambda)}
\newcommand\SVl{S_V(\lambda)}
\newcommand\Sjr{S_j^{\reduced}}
\newcommand\Rkp{{\mathcal R}^k_+}
\newcommand\Rkm{{\mathcal R}^k_-}
\newcommand\Rkpm{{\mathcal R}^k_{\pm}}
\newcommand\Phys{{\mathcal P}}
\newcommand\Pc{\overline{\mathcal P}}
\newcommand\pip{\pi^{\perp}}
\newcommand\pipa{\pi_\partial}
\newcommand\gammapa{\gamma_\partial}
\newcommand\pipah{\hat\pi_\partial}
\newcommand\pit{\tilde\pi}
\newcommand\xit{\tilde\xi}
\newcommand\zetat{\tilde\zeta}
\newcommand\etat{\tilde\eta}
\newcommand\sigmat{\tilde\sigma}
\newcommand\sigmahat{\hat\sigma}
\newcommand\thetat{\tilde\theta}
\newcommand\psit{\tilde\psi}
\newcommand\phit{\tilde\phi}
\newcommand\chit{\tilde\chi}
\newcommand\rhot{\tilde\rho}
\newcommand\xib{\bar\xi}
\newcommand\zetab{\bar\zeta}
\newcommand\thetab{\bar\theta}
\newcommand\etab{\bar\eta}
\newcommand\iotal{\iota_{\lambda}}
\newcommand\rhoat{\rhot_{\alpha_1}}
\newcommand\Lambdat{\tilde\Lambda}
\newcommand\Lambdati{\tilde\Lambda^{\text{in}}}
\newcommand\Lambdato{\tilde\Lambda^{\text{out}}}
\newcommand\Lambdatp{\tilde\Lambda^{\text{prop}}}
\newcommand\Lambdai{\Lambda^{\text{in}}}
\newcommand\Lambdao{\Lambda^{\text{out}}}
\newcommand\poles{\Lambda'}
\newcommand\rpoles{\Lambda_p}
\newcommand\thresholds{\Lambda}
\newcommand\Vt{\tilde V}
\newcommand\It{\tilde I}
\newcommand\half{{\frac{1}{2}}}
\newcommand\sigmah{\sigma^{1/2}}
\newcommand\bX{\partial X}
\newcommand\bXb{\partial \Xb}
\newcommand\Deltabt{\tilde\Delta_0}
\newcommand\strip{\Omega_T}
\newcommand\Kf{K^{\flat}}
\newcommand\Gs{G^{\sharp}}
\newcommand\Gt{\tilde G}
\newcommand\Osc{\sci\Omega}
\newcommand\OSc{{}^\Scl\Omega}
\newcommand\Osch{\sci\Omega^{\half}}
\newcommand\Oscmh{\sci\Omega^{-\half}}
\newcommand\Isc{I_{sc}}
\newcommand\os{{\text{os}}}
\newcommand\Qzl{Q^0_{-\lambda}}
\newcommand\Lie{{\mathcal L}}
\newcommand\bl{{\text b}}
\newcommand\scl{{\text{sc}}}
\newcommand\sccl{{\text{scc}}}
\newcommand\Scl{{\text{Sc}}}
\newcommand\ScLl{{\text{Sc,L}}}
\newcommand\ScRl{{\text{Sc,R}}}
\newcommand\Sccl{{\text{Scc}}}
\newcommand\sus{{\text{sus}}}
\newcommand\ssl{{\text{ee}}}
\newcommand\bzl{{\text{b0}}}
\newcommand\XXb{X^2_\bl}
\newcommand\XXbt{\Xt^2_\bl}
\newcommand\XXsc{X^2_\scl}
\newcommand\XXsct{\Xt^2_\scl}
\newcommand\XXSc{X^2_\Scl}
\newcommand\XXSct{\Xt^2_\Scl}
\newcommand\XXScL{X^2_\ScLl}
\newcommand\XXScR{X^2_\ScRl}
\newcommand\MMsc{M^2_\scl}
\newcommand\Deltab{\Delta_\bl}
\newcommand\Deltasc{\Delta_\scl}
\newcommand\DeltaSc{\Delta_\Scl}
\newcommand\DeltaScL{\Delta_\ScLl}
\newcommand\DeltaScR{\Delta_\ScRl}
\newcommand\prs{\sigma}
\newcommand\Nsc{N_\scl}
\newcommand\Nscp{N_{\scl,p}}
\newcommand\Nff{N_{\ff}}
\newcommand\Nffz{N_{\ff,0}}
\newcommand\Nffzp{N_{\ff,0,p}}
\newcommand\Nffl{N_{\ff,l}}
\newcommand\Nffml{N_{\ff,-l}}
\newcommand\Nmf{N_{\mf}}
\newcommand\Nmfz{N_{\mf,0}}
\newcommand\Nmfl{N_{\mf,l}}
\newcommand\Nmfml{N_{\mf,-l}}
\newcommand\ffb{\operatorname{bf}}
\newcommand\Ffb{\operatorname{bf'}}
\newcommand\ffsc{\operatorname{sf}}
\newcommand\ffSc{\operatorname{sf_C}}
\newcommand\Ffsc{\operatorname{sf'}}
\newcommand\rff{\rho_{\ff}}
\newcommand\rmf{\rho_{\mf}}
\newcommand\rffb{\rho_{\ffb}}
\newcommand\rffsc{\rho_{\ffsc}}
\newcommand\rFfsc{\rho_{\Ffsc}}
\newcommand\rffSc{\rho_{\ffSc}}
\newcommand\rinf{\rho_{\infty}}
\newcommand\CL{C_L}
\newcommand\CR{C_R}
\newcommand\betab{\beta_\bl}
\newcommand\betasc{\beta_\scl}
\newcommand\betaSc{\beta_\Scl}
\newcommand\BetaSc{\bar\beta_\Scl}
\newcommand\betaScL{\beta_\ScLl}
\newcommand\betaScR{\beta_\ScRl}
\newcommand\ScT{{}^\Scl T^*}
\newcommand\SccT{{}^\Scl \bar T^*}
\newcommand\ScS{{}^\Scl S^*}
\newcommand\Tb{{}^\bl T}
\newcommand\Tss{{}^\ssl T}
\newcommand\Tsc{{}^\scl T}
\newcommand\TSc{{}^\Scl T}
\newcommand\CSc{C_\Scl}
\newcommand\Lambdasc{{}^\scl\Lambda}
\newcommand\XXXb{X^3_\bl}
\newcommand\XXXsc{X^3_\scl}
\newcommand\XXXSc{X^3_\Scl}
\newcommand\XXXScO{X^3_{\Scl,O}}
\newcommand\XXXScF{X^3_{\Scl,F}}
\newcommand\XXXScS{X^3_{\Scl,S}}
\newcommand\XXXScC{X^3_{\Scl,C}}
\newcommand\KDsc{\operatorname{KD^{\half}_\scl}}
\newcommand\KDSc{\operatorname{KD^{\half}_\Scl}}
\newcommand\KDScEF{\operatorname{KD^{E,F}_\Scl}}
\newcommand\Oh{\operatorname{\Omega^{\half}}}
\newcommand\WFSc{\WF_\Scl}
\newcommand\WFtSc{\WF_{\text 3sc}}
\newcommand\WFScmf{\WF_{\Scl,\mf}}
\newcommand\WFScff{\WF_{\Scl,\ff}}
\newcommand\WFScs{\WF_{\Scl,\prs}}
\newcommand\WFScp{\WF'_\Scl}
\newcommand\WFScmfp{\WF'_{\Scl,\mf}}
\newcommand\WFScffp{\WF'_{\Scl,\ff}}
\newcommand\WFScsp{\WF'_{\Scl,\prs}}
\newcommand\Diffscc{\Diff_\sccl}
\newcommand\DiffSc{\Diff_\Scl}
\newcommand\Diffss{\Diff_\ssl}
\newcommand\DiffScc{\Diff_\Sccl}
\newcommand\DiffscI{\Diff_{\scl,\text{I}}}
\newcommand\VscI{\Vf_{\scl,\text{I}}}
\newcommand\DiffsV{\operatorname{Diff}_{\sus(V)}}
\newcommand\DiffsVsc{\operatorname{Diff}_{\sus(V),\scl}}
\newcommand\DiffsVCsc{\operatorname{Diff}_{\sus(V)-C,\scl}}   
\newcommand\Psisc{\Psop_\scl}
\newcommand\Psiscc{\Psop_\sccl}
\newcommand\Psiss{\Psop_\ssl}
\newcommand\Psisch{\Psop_{\scl,h}}
\newcommand\Psiscch{\Psop_{\sccl,h}}
\newcommand\PsiSc{\Psop_\Scl}
\newcommand\PsiScph{\Psop_{\Scl,\phi}}
\newcommand\PsiScra{\Psop_{\Scl,\rho^\sharp_a}}
\newcommand\PsiScc{\Psop_\Sccl}
\newcommand\PsiSccml{\Psop^{m,l}_\Sccl}
\newcommand\PsiScxx{\Psop^{*,*}_\Scl}
\newcommand\PsiScml{\Psop^{m,l}_\Scl}
\newcommand\PsiScmz{\Psop^{m,0}_\Scl}
\newcommand\PsiScmmz{\Psop^{-m,0}_\Scl}
\newcommand\PsiSckz{\Psop^{k,0}_\Scl}
\newcommand\PsiScmmml{\Psop^{-m,-l}_\Scl}
\newcommand\Psiscmkk{\Psop^{-k,k}_\scl}
\newcommand\Psiscmmmkk{\Psop^{-m-k,k}_\scl}
\newcommand\Psiscmoo{\Psop^{-1,1}_\scl}
\newcommand\Psiscmz{\Psop^{m,0}_\scl}
\newcommand\Psiscmmz{\Psop^{-m,0}_\scl}
\newcommand\PsiSckmkl{\Psop^{km,kl}_\Scl}
\newcommand\PsiScmplp{\Psop^{m',l'}_\Scl}
\newcommand\PsiScmmpllp{\Psop^{m+m',l+l'}_\Scl}
\newcommand\Psiscml{\Psop^{m,l}_\scl}
\newcommand\PsiScid{\Psop^{0,0}_\Scl}
\newcommand\PsiSczo{\Psop^{0,1}_\Scl}
\newcommand\PsiScmii{\Psop^{-\infty,\infty}_\Scl}
\newcommand\PsiScmiz{\Psop^{-\infty,0}_\Scl}
\newcommand\PsiScmoo{\Psop^{-1,1}_\Scl}
\newcommand\PsisCid{\Psop^{0,0}_{\scl-C}}
\newcommand\PsisC{\Psop_{\scl-C}}
\newcommand\Psiinf{\Psop_{\infty}}
\newcommand\Psiinfid{\Psop_{\infty}^0}
\newcommand\PsiFinf{\Psop_{\infty-\Fr}}
\newcommand\PsisVscml{\Psop^{m,l}_{\sus(V),\scl}}
\newcommand\PsisVsc{\Psop_{\sus(V),\scl}}
\newcommand\PsisVpsc{\Psop_{\sus(V_p),\scl}}
\newcommand\PsisVCSc{\Psop_{\sus(V)-C,\scl}}
\newcommand\SFinf{S_{\infty-\Fr}}
\newcommand\YsVC{Y^2_{\sus(V)-C,\scl}}
\newcommand\ffYsc{\ffsc_{\sus(V)}}
\newcommand\SXC{S(X;C)}
\newcommand\Ios{I_{\text{os}}}
\newcommand\pbL{\pi^2_{\bl,{\text L}}}
\newcommand\pbR{\pi^2_{\bl,{\text R}}}
\newcommand\pscL{\pi^2_{\scl,{\text L}}}
\newcommand\pscR{\pi^2_{\scl,{\text R}}}
\newcommand\PbO{\pi^3_{\bl,{\text O}}}
\newcommand\PscO{\pi^3_{\scl,{\text O}}}
\newcommand\PScO{\pi^3_{\Scl,{\text O}}}
\newcommand\PScF{\pi^3_{\Scl,{\text F}}}
\newcommand\PScC{\pi^3_{\Scl,{\text C}}}
\newcommand\PScS{\pi^3_{\Scl,{\text S}}}
\newcommand\pScL{\pi^2_{\Scl,{\text L}}}
\newcommand\pScR{\pi^2_{\Scl,{\text R}}}
\newcommand\CLF{\CL^F}
\newcommand\CLO{\CL^O}
\newcommand\CLS{\CL^S}
\newcommand\CLC{\CL^C}
\newcommand\DeltaYb{\Delta_{\bl,Y}}
\newcommand\DeltaYsc{\Delta_{\sus-\scl}}
\newcommand\diag{\operatorname{diag}}
\newcommand\Diag{\operatorname{Diag}}
\newcommand\Vf{{\mathcal V}}
\newcommand\Vb{{\mathcal V}_{\bl}}
\newcommand\Vsc{{\mathcal V}_{\scl}}
\newcommand\Vss{{\mathcal V}_{\ssl}}
\newcommand\VSc{{\mathcal V}_{\Scl}}
\newcommand\VfI{\Vf_{\text{I}}}
\newcommand\VfIq{\Vf_{\text{I},q}}
\newcommand\scH{{}^\scl H}
\newcommand\scHg{\scH_g}
\newcommand\Hss{H_\ssl}
\newcommand\Hbz{H_\bzl}
\newcommand\xh{\hat x}
\newcommand\Yh{\hat Y}
\newcommand\Zh{\hat Z}
\newcommand\Yb{\bar Y}
\newcommand\hb{\bar h}
\newcommand\xih{\hat\xi}
\newcommand\etah{\hat\eta}
\newcommand\muh{\hat\mu}
\newcommand\mub{\bar\mu}
\newcommand\nub{\bar\nu}
\newcommand\mubh{\widehat{\bar\mu}}
\newcommand\yb{\bar y}
\newcommand\zb{\bar z}
\newcommand\ub{\bar u}
\newcommand\Qb{\bar Q}
\newcommand\Wbp{{\bar W}^\perp}
\newcommand\Wp{W^\perp}
\newcommand\Kt{\tilde K}
\newcommand\Wt{\tilde W}
\newcommand\Ut{\tilde U}
\newcommand\yt{\tilde y}
\newcommand\ut{\tilde u}
\newcommand\vt{\tilde v}
\newcommand\ft{\tilde f}
\newcommand\htil{\tilde h}
\newcommand\St{\tilde S}
\newcommand\Pt{\tilde P}
\newcommand\Rt{\tilde R}
\newcommand\qt{\tilde q}
\newcommand\Qt{\tilde Q}
\newcommand\Xb{\overline{X}}
\newcommand\lambdat{\tilde\lambda}
\newcommand\betat{\tilde\beta}
\newcommand\epst{\tilde\epsilon}
\newcommand\ep{\epsilon}
\newcommand\bt{\tilde b}
\newcommand\Xt{\widetilde X}
\newcommand\Mt{\widetilde M}
\newcommand\At{\tilde A}
\newcommand\Et{\tilde E}
\newcommand\Ht{\tilde H}
\newcommand\at{\tilde a}
\newcommand\Ct{\tilde C}
\newcommand\pih{\hat\pi}
\newcommand\Rh{\hat R}
\newcommand\Ah{\hat A}
\newcommand\Bh{\hat B}
\newcommand\Ch{\hat C}
\newcommand\Gh{\hat G}
\newcommand\Hh{\hat H}
\newcommand\Qh{\hat Q}
\newcommand\Ph{\hat P}
\newcommand\Nh{\hat N}
\newcommand\Sh{\hat S}
\newcommand\Gcal{{\mathcal G}}
\newcommand\GcalC{{\mathcal G}_C}
\newcommand\Jcal{{\mathcal J}}
\newcommand\JcalC{{\mathcal J}_C}
\setcounter{secnumdepth}{3}
\newtheorem{lemma}{Lemma}[section]
\newtheorem{prop}[lemma]{Proposition}
\newtheorem{thm}[lemma]{Theorem}
\newtheorem{cor}[lemma]{Corollary}
\newtheorem{result}[lemma]{Result}
\newtheorem*{thm*}{Theorem}
\newtheorem*{prop*}{Proposition}
\newtheorem*{conj*}{Conjecture}
\numberwithin{equation}{section}
\theoremstyle{remark}
\newtheorem{rem}[lemma]{Remark}
\theoremstyle{definition}
\newtheorem{Def}[lemma]{Definition}
\newtheorem*{Def*}{Definition}
\def\signature#1#2{\par\noindent#1\dotfill\null\\*
{\raggedleft #2\par}}